\documentclass[aos,MSNbibl,seceqn,citesort,dvips]{arximspdf}

\doi{10.1214/11-AOS963} 
\volume{40}
\issue{1}
\pubyear{2012}
\firstpage{494}
\lastpage{528}

\makeatletter

\newcommand{\xrightarrow}[1]{\stackrel{#1}{\rightarrow}}

\newtheorem{theorem}{Theorem}[section]
\newtheorem{cor}[theorem]{Corollary}
\newtheorem{lem}[theorem]{Lemma}
\newtheorem{prop}[theorem]{Proposition}

\newproclaim{ass}[theorem]{Assumption}
\newproclaim{rem}[theorem]{Remark}

\newcommand{\argmin}{\mathop{\arg\min}}

\newcommand{\xiv}{\bolds\xi}
\newcommand{\etav}{\bolds\eta}
\newcommand{\ind}{\mathbf{1}}
\newcommand{\R}{\mathbb{R}}
\newcommand{\Z}{\mathbb{Z}}
\newcommand{\Cov}{\operatorname{Cov}}
\newcommand{\Var}{\operatorname{Var}}
\newcommand{\E}{\mathbb{E}}
\newcommand{\N}{\mathbb{N}}
\newcommand{\F}{\mathcal{F}}
\newcommand{\OO}{\mathcal{O}}
\newcommand{\OOO}{\mbox{\scriptsize$\mathcal{O}$}}
\newcommand{\Gaussian}{\mathcal{N}}

\makeatother

\begin{document}
\begin{frontmatter}

\title{Simultaneous confidence bands for Yule--Walker estimators and order selection}
\runtitle{Confidence bands and order selection}

\begin{aug}
\author[A]{\fnms{Moritz} \snm{Jirak}\corref{}\ead[label=e1]{m0ritz@yahoo.com}}
\runauthor{M. Jirak}
\affiliation{Graz University of Technology}
\address[A]{Institute of Statistics\\
Graz University of Technology\\
M\"unzgrabenstra{\ss}e 11, A-8010 Graz\\
Austria\\
\printead{e1}} 
\end{aug}

\received{\smonth{12} \syear{2010}}
\revised{\smonth{9} \syear{2011}}

%
\begin{abstract}
Let $\{X_k, k\in{\mathbb Z}\}$ be an autoregressive process of order
$q$. Various estimators for the order $q$ and the parameters ${\bolds
\Theta}_q = (\theta_1,\ldots,\theta_q)^T$ are known; the order is usually
determined with Akaike's criterion or related modifications, whereas
Yule--Walker, Burger or maximum likelihood estimators are used for the
parameters ${\bolds\Theta}_q$. In this paper, we establish simultaneous
confidence bands for the Yule--Walker estimators
$\widehat{\theta}_{i}$; more precisely, it is shown that the limiting
distribution of ${\max_{1 \leq i \leq d_n}}|\widehat{\theta
}_{i} -
\theta_{i}|$ is the Gumbel-type distribution $e^{-e^{-z}}$, where
$q \in\{0,\ldots,d_n\}$ and $d_n = \OO(n^{\delta} )$,
$\delta
> 0$. This allows to modify some of the currently used criteria (AIC,
BIC, HQC, SIC), but also yields a~new class of consistent estimators
for the order $q$. These estimators seem to have some potential, since
they outperform most of the previously mentioned criteria in a~small
simulation study. In particular, if some of the parameters
$\{\theta_i\}_{1 \leq i \leq d_n}$ are zero or close to zero,
a~significant improvement can be observed. As a~byproduct, it is shown
that BIC, HQC and SIC are consistent for $q \in\{0,\ldots,d_n\}$ where
$d_n = \OO(n^{\delta} )$.
\end{abstract}

%

\end{frontmatter}

\section{Introduction}\label{secintro}

$\!\!\!$Let $\{X_k\}_{k \in\Z}$ be a~$q$th-order autoregressive process
$\operatorname{AR}(q)$
with coefficient vector ${\bolds\Theta}_q \in\R^{q}$. A~considerable
literature in the past years dealt with various aspects and problems on
$\operatorname{AR}(q)$-processes; see, for instance,
\cite{andersentimeseries,timeseriesbrockwell,hannantimeseries,luetkepohl}
and the references therein. More recently, people have moved on to more
complicated models such as ARCH~\cite{bollerslevarch,engle},
GARCH~\cite{bollerslevgarch} and related models, which again have been
extended in many different directions. However, in many applications,
$\operatorname{AR}(q)$-processes still form the backbone and are often used as first
approximations for further analysis; in particular, many estimation and
fitting procedures can be based on preliminary $\operatorname{AR}(0q)$ approximations.
This includes, for instance, ARMA, ARCH and GARCH models
\cite{bhansali1991,gourieuxbook,hannanarma}. Thus, $\operatorname{AR}(q)$
processes have moved from the spotlight to the backstage area, yet
their significance remains unchallenged.

When fitting an $\operatorname{AR}(q)$ model, two important questions arise: how to
choose the order $q$, and having done so, which estimators are to be
used. Naturally, these two problems can hardly be separated and are
often dealt with simultaneously, or at least so in preliminary
estimates. An extensive literature has evolved around these two issues.
Pioneering contributions in this direction are due to Akaike
\mbox{\cite{akaike1969,akaike1977}}, Mallows
\cite{mallows1964,mallows1966}, Walker~\cite{Walker1931} and
Yule~\cite{yule1927}; for more details we refer to
\cite
{andersentimeseries,boxbook,timeseriesbrockwell,hannantimeseries,luetkepohl}
and the references there. In order to be able to describe some of the
basic results, we recall that an $\operatorname{AR}(q)$ process $\{X_k\}_{k \in\Z}$ is
defined through the recurrence relation
%
%
\begin{equation}\label{eqeqarp}
X_k = \theta_{1} X_{k-1} + \cdots+ \theta_{q} X_{k-q} + \varepsilon_{k},
\end{equation}
where it is often assumed that $\{\varepsilon_k\}_{k \in\Z }$ is a~mean-zero i.i.d. sequence. Let $\phi_{h} = \E(X_k X_{k + h})$, $k,h
\in\Z$, be the covariance function. A~natural estimate for $\phi_h$ is
the sample covariance $\widehat{\phi}_{n,h} = \frac{1}{n} \sum_{i = h +
1}^n X_i X_{i - h}$. Depending on the magnitude of $h$, a~different
normalization, such as $(n - h)^{-1}$, is sometimes more convenient.
Denote with ${\bolds \Theta}_{{q}} = (\theta_1,\ldots,\theta_{q})^T$
the parameter vector and put ${\bolds\Phi}_{q} =
(\phi_{1},\ldots,\phi_{q} )^T$, and let ${\bolds\Gamma}_{q} =
(\phi_{|i-j|} )_{1 \leq i,j \leq q}$ be the ($q\times q$)-dimensional
covariance matrix. Then it follows from~(\ref{eqeqarp}) that $
{\bolds\Gamma}_{q} {\bolds\Theta}_{{q}} = {\bolds\Phi}_{q}; $ hence
a\vspace*{1pt} natural idea is to replace the corresponding quantities
by estimators $\widehat{\bolds\Phi}_{q} = (\widehat{\phi
}_{n,1},\ldots,\widehat{\phi}_{n,{q}} )^T$, $\widehat{\bolds\Gamma}_{q}
= (\widehat{\phi}_{n,|i-j|} )_{1 \leq i,j \leq{q}}$, and thus define
the estimator $\widehat{\bolds\Theta}_{{q}} = (\widehat{\theta
}_1,\ldots,\widehat {\theta}_{q})^T$ via
%
%
\begin{equation}\label{defnsigma}
\widehat{\bolds\Gamma}_{q}^{-1} \widehat{\bolds\Phi}_q = \widehat{\bolds
\Theta
}_q \quad\mbox{and}\quad \widehat{\sigma}^2(q) = \widehat{\phi}_0 -
\widehat
{\bolds\Theta}_{{q}}^T \widehat{\bolds\Phi}_{q},
\end{equation}
where $\sigma^2 = \E(\varepsilon_0^2)$. These estimators are
commonly referred to as the Yule--Walker estimators, and they have some
remarkable properties. For example, if $\{X_k\}_{k \in\Z}$
is causal, then the fitted model
\[
X_k = \widehat{\theta}_{1} X_{k-1} + \cdots+ \widehat{\theta}_{p}
X_{k-q} + \varepsilon_{k}
\]
is still causal; see, for instance,~\cite{timeseriesbrockwell} and
\cite{nakaregion}. Another interesting feature is that even though
the Yule--Walker estimators are obtained via moment matching methods,
their variance is asymptotically equivalent with those obtained via a~maximum likelihood approach.
More precisely, for $m \geq q$ it holds that
%
%
\begin{equation}\label{eqcltholds}
\sqrt{n} (\widehat{\bolds\Theta}_{{m}} - {\bolds\Theta
}_{{m}})
\xrightarrow{d} \Gaussian(0, \sigma^2 {\bolds\Gamma
}_m^{-1}),
\end{equation}
where ${\bolds\Theta}_{{m}} = (\theta_1,\ldots,\theta_{q},
0,\ldots,0
)^T$; see, for instance,~\cite{timeseriesbrockwell}. These asymptotic
results form the basis for earlier estimation methods of the order $q$
\cite{Quenouille1947,walker1952,whittle1952}, which focused on a~fixed, finite number of possible orders and consist of
multiple-testing-procedures, which in practice leads to the difficulty
of having a~required level. On the other hand, as it was pointed out by
Shwarz~\cite{Shwarz1978}, a~direct likelihood approach fails, since
it invariably chooses the highest possible dimension. Akaike
\cite{akaike1969} and Mallows~\cite{mallows1964,mallows1966}, developed
a~different approach, which is based on a~``generalized likelihood
function.'' Shibata~\cite{shibata} investigated the asymptotic
distribution and showed that the estimator based on~(\ref{defnaik}) is
not consistent. This issue was successfully dealt with by
Akaike~\cite{akaike1977} (BIC), Hannan and Quinn~\cite{hannanarp} (HQC),
Parzen~\cite{parzen1974}, Rissanen~\cite{Rissanen1978} and
Shwarz~\cite{Shwarz1978} (SIC), who introduced consistent modifications (Parzen's
CAT-criterion is conceptually different). For more recent advances and
generalizations, see, for instance, Barron et al.~\cite{barron1999},
Foster and George~\cite{foster1994}, Shao~\cite{shao1997} and the
detailed review on model selection given by Leeb and P\"{o}tscher
\cite{leep2008}. A~particularly interesting direction addresses
$\operatorname{AR}(\infty
)$ approximations; recent contributions are due to Bickel and Yel
\cite{bickel2011} and Ing and Wei \mbox{\cite{ingwei2003,ingwei2005}}. Here
and now, we will content ourselves with briefly discussing Akaike's
approach and closely related criteria. Akaike's generalized likelihood
function leads to the expression
%
%
\begin{equation}\label{defnaik}
\mathrm{AIC}(m) = n \log\widehat{\sigma}^2(m) + 2 m,
\end{equation}
where $n$ is the sample size and $\widehat{\sigma}^2(m)$ is as in
(\ref{defnsigma}). An estimator for the order $q$ is then obtained by
minimizing $\mathrm{AIC}(m)$, $m \in\{0,1,\ldots,K\}$, for some
predefined $0 \leq q \leq K$. Consistent modifications are obtained by
inserting an increasing sequence~$C_n$, and $\mathrm{AIC}(m)$ then becomes
%
%
\begin{equation}
\widetilde{\mathrm{AIC}}(m) = n \log\widehat{\sigma}^2(m) + 2 C_n
m,\qquad
m \in\{0,1,\ldots,K\}.
\end{equation}
Most modifications result in $C_n = \OO( \log n )$, even
though the arguments are sometimes quite different. A~notable exception
is the idea of Hannan and Quinn~\cite{hannanarp}, who successfully
employed the LIL to obtain $C_n = \OO(\log\log n )$.

The aim of this paper is to introduce a~different approach, based on
the quantity ${\max_{1 \leq i \leq d_n}} |\widehat{\theta}_{{i}} -
{\theta}_{{i}}|$, where $d_n$ is an increasing function in $n$. It
is shown, for instance, that, appropriately normalized, this expression
converges weakly to a~Gumbel-type distribution. On one hand, this
allows to construct simultaneous confidence bands for the Yule--Walker
estimators~$\widehat{\bolds\Theta}_{{d_n}}$, but also permits us to
construct a~variety of different, consistent estimators for the order
$q$ of an autoregressive process. The asymptotic distribution of such a~particular estimator is also derived. As a~byproduct, it is shown that
known consistent criteria such as BIC, SIC and HQC are also consistent
if the parameter space is increasing; that is, consistency even holds
if $q \in\{0,\ldots,d_n\}$, where $d_n = \OO(n^{\delta} )$.
This partially gives answers to questions raised by Hannan and Quinn
\cite{hannanarp} and Shibata~\cite{shibata}, and extends results
given by An et al.~\cite{hannan1982}.
In addition, the general method seems to be very useful for model
fitting for subset autoregressive processes (see, e.g.,~\cite{mcleod2008}),
which is highlighted in Remark~\ref{remautosubset}
and Section~\ref{secnumerical}. A~more thorough treatment of this
issue is postponed to a~subsequent paper.

\section{Main results}\label{secmain}

We will frequently use the following notation. For a~vector $x =
(x_1,\ldots,x_d)^T$, we put $\|x\|_{\infty} = {\max_{1 \leq i \leq
d}}|x_i|$, and for a~matrix ${\mathbf A} = (a_{i,j})_{\{
{1 \leq i \leq r, 1 \leq j \leq s}\}}$, $r,s \in\N$
we denote with
%
%
\begin{equation}
\|{\mathbf A}\|_{\infty} = \max\{ {\mathbf A} x | x \in
\R
^s, \|x\|_{\infty} = 1 \} = \max_{1 \leq i \leq r} \sum_{j =
1}^s |a_{i,j}|
\end{equation}
the usual induced matrix norm. In addition, we will use the
abbreviation $\|\cdot\|_p =(\E(|\cdot|^p))^{1/p}$, $p < \infty$. The
main results involve an array of $\operatorname{AR}(q)$ processes; more precisely, we
consider the family of $\operatorname{AR}(d_n)$ processes $\{X_k^{(r)}\}_{k
\in\Z}$, $1 \leq r \leq d_n$, where $d_n = \OO(n^{\delta}
)$ (more details are given later). Since we are always only dealing
with a~single member of this array, the index $(r)$ is dropped for
convenience, and we just consider an $\operatorname{AR}(d_n)$ process $\{
X_k
\}_{ k \in\Z}$, keeping in mind that the parameters $\{\theta
_i\}_{1 \leq i \leq d_n}$ may depend on $n$. This implies that
$X_k$ satisfies the recurrence relation
%
%
\begin{equation}
X_k = \theta_{1} X_{k-1} + \cdots+ \theta_{d_n} X_{k-d_n} + \varepsilon
_{k},\qquad k \in\Z,
\end{equation}
where $\{\varepsilon_k\}_{ k \in\Z}$ defines the usual
innovations. Note that $d_n$ does not need to reflect the actual order
$q$ of the $\operatorname{AR}(d_n)$ process, as we do not require that $\{\theta
_i\}_{1 \leq i \leq d_n}$ are all different from zero. All of the
results are derived under the following assumption regarding the
$\operatorname{AR}(d_n)$ process $\{X_k\}_{ k \in\Z}$.
%
%
\begin{ass}\label{assmain}
$\{X_k\}_{ k \in\Z}$ admits a~causal representation $X_k =\break \sum_{i =
0}^{\infty} \alpha_i \times\varepsilon_{k - i}$, such that:
\begin{itemize}
\item$\sup_n \Psi(m) = \OO(m^{-\vartheta})$,
$\vartheta>
0$, where $\Psi(m) := {\sum_{i = m}^{\infty}} |\alpha_i|$,
\item$\{\varepsilon_k\}_{k \in\Z}$ is a~mean-zero i.i.d.
sequence of random variables, such that $\|\varepsilon_k\|_p <
\infty$ for some $p > 4$, $\|\varepsilon_k\|_2^2 = \sigma^2 >
0$, $k \in\Z$,
\item${\sup_n \sum_{i = 1}^{\infty} }|\theta_i| < \infty$,
$|\theta_n| = \OOO
((\log n)^{-1})$.
\end{itemize}
\end{ass}

In accordance with the previously established notation, we introduce
the inverse and estimated inverse matrix
%
%
\begin{equation}
{\bolds\Gamma}_{d_n}^{-1} = (\gamma_{i,j}^* )_{1 \leq i,j
\leq
{d_n}},\qquad
\widehat{\bolds\Gamma}_{d_n}^{-1} = (\widehat{\gamma}_{i,j}^*
)_{1 \leq i,j \leq{d_n}}.
\end{equation}
In addition, we will use the convention
that $\theta_0 = \widehat{\theta}_0 = -1$. We can now formulate our
main result.
%
%
\begin{theorem}\label{thmyeconvrgextreme}
Let $\{X_k\}_{k \in\Z}$ be an $\operatorname{AR}(d_n)$ process satisfying
Assumption~\ref{assmain}.
Suppose that $d_n \to\infty$ as $n$ increases, with $d_n = \OO
(n^{\delta})$ such that
%
%
\begin{equation}
0 < \delta< \min\{1/2, \vartheta p/2\},\qquad (1 - 2
\vartheta) \delta< (p - 4)/p.
\end{equation}
If we have in addition that ${\inf_{h}}|\gamma_{h,h}^* | > 0$,
then for $z \in\R$
\[
P\Bigl(a_n^{-1}\Bigl({\sqrt{n}\max_{1 \leq i \leq d_n}} |
(\widehat{\gamma}_{i,i}^* \widehat{\sigma}^2(d_n))^{-1/2}(\widehat
{\theta}_{i} - \theta_{i})| - b_{n}\Bigr) \leq z \Bigr)
\rightarrow\exp(- e^{-z}),
\]
where $a_{n}= (2 \log d_n)^{-1/2}$ and $b_{n} = (2 \log d_n)^{1/2} - (8
\log d_n)^{-1/2} (\log\log d_n + 4 \pi- 4)$.
\end{theorem}
%
%
\begin{rem}\label{remdiscusscondis}
Condition $\inf_{h}|\gamma_{h,h}^* | > 0$ may be explicitly
expressed in terms of $\{\theta_i\}_{1 \leq i \leq d_n}$ [see
(\ref{eqgammainversematrixexact})], and is quite general. In fact,
it is only needed to control or exclude possible pathological cases.
\end{rem}
%
%
\begin{rem}
Note that if we have $|\alpha_i| = \OO(i^{-3/2})$, then
$\vartheta\geq1/2$. Hence condition $p > 4$ implies that we may
choose $\delta$ arbitrarily close to $1/2$, which essentially results
in $d_n = \OOO(\sqrt{n})$.
\end{rem}

The above remark indicates that we may obtain simple bounds for $d_n$,
provided that we can control $\alpha_i$ asymptotically. If the
cardinality of the set $\{1 \leq i \leq d_n |\theta_i \neq0
\}$ tends to infinity as $n$ increases, then establishing general
and simple conditions on the relation between $\{\theta_i\}
_{1 \leq i \leq d_n}$ and $\Psi(m)$ seems to be very difficult. One
may, however, obtain the following corollary.
%
%
\begin{cor}\label{corsimple}
Suppose that $\{\varepsilon_k\}_{k \in\Z}$ is a~mean-zero
i.i.d. sequence of random variables, such that $\|\varepsilon
_k\|
_p < \infty$ for some $p > 4$, $\|\varepsilon_k\|_2^2 =
\sigma^2
> 0$, $k \in\Z$, and that one of the following conditions holds:
\begin{longlist}
\item ${\sup_n \sum_{i = 1}^{d_n}} |\theta_i | < 1$,
\item $\theta_i = 0$, $q < i \leq d_n$ for some fixed $q \in\N$
which does not depend on $n$.
\end{longlist}
Then the conditions of Theorem~\ref{thmyeconvrgextreme} are
satisfied and we can choose any $d_n = \OO(n^{\delta})$ with
$\delta< 1/2$.
\end{cor}
%
%
\begin{rem}\label{rempractice}
The rate of convergence to an extreme-value type distribution as given
in Theorem~\ref{thmyeconvrgextreme} can be rather slow; see, for
instance,~\cite{balkema1990,omey1989}. Hence, in view of~(\ref{eqcltholds}) (and Theorem~\ref{thmgaussianapprox}), it may be more
appropriate to use the approximation
\[
P\Bigl( {\max_{1 \leq i \leq n}}\bigl|\sqrt{n}(\widehat{\gamma}_{i,i}^*
\widehat{\sigma}^2)^{-1/2}(\widehat{\theta}_{i} - \theta_{i})\bigr|
\leq x \Bigr) \approx P(\|{\bolds\xi_{d_n}}\|
_{\infty}
\leq x )
\]
in practice, where $\xiv_{d_n} = (\xi_{n,1},\ldots, \xi_{n,d_n}
)^T$ is a~$d_n$-dimensional mean-zero Gaussian random vector with the
same covariance structure. Corresponding quantiles can be obtained, for
instance, via a~Monte Carlo technique. However, if $d_n$ is
sufficiently large, one has that
\[
P(\|\xiv_{d_n}\|_{\infty} \leq x ) \approx
P
(\| \etav_{d_n}\|_{\infty} \leq x ),
\]
where $\etav_{d_n} = (\eta_{n,1},\ldots, \eta_{n,d_n} )^T$
is a~sequence of i.i.d. mean-zero Gaussian random variables with unit
variance. A~bound for the error can be given by using the techniques
developed by Berman~\cite{Berman1964} and Deo~\cite{deoabs}; see
also the proof of Theorem~\ref{thmestdistrib}.
\end{rem}

The above results allow us to construct the simultaneous confidence bands
%
%
\begin{eqnarray}\label{confidencemax}
\mathcal{M}_1(d_n) &=& \Bigl\{{\bolds\Theta}_{{d_n}} \in\R^{d_n}
\big|\nonumber\\[-8pt]\\[-8pt]
&&\hspace*{7pt}a_n^{-1}\Bigl(\sqrt{n}{\max_{1 \leq i \leq d_n}} |
(\widehat{\gamma}_{i,i}^*)^{-1/2}
(\widehat{\theta}_{i} - \theta
_{i})| - b_{n}\Bigr) \leq  \sqrt{\widehat{\sigma}^2(d_n)}
V_{1 - \alpha}\Bigr\},\hspace*{-30pt}\nonumber
\end{eqnarray}
where $V_{1 - \alpha}$ denotes the $1 - \alpha$ quantile of the
Gumbel-type distribution given above. In the literature
\cite{andersentimeseries,timeseriesbrockwell,hannantimeseries} one
often finds the confidence ellipsoids
%
%
\begin{eqnarray}\label{confidenceellips}
\mathcal{M}_2(m) &=& \{{\bolds\Theta}_{{m}} \in\R^{m}
|\nonumber\\[-8pt]\\[-8pt]
&&\hspace*{4pt}(\widehat{\bolds\Theta}_{{m}} - {\bolds\Theta}_{{m}})
\widehat
{\bolds\Gamma}_m (\widehat{\bolds\Theta}_{{m}} - {\bolds\Theta
}_{{m}})^T \leq n^{-1} \widehat{\sigma}^2(m) \chi_{1 - \alpha
}^2(m)\},\nonumber
\end{eqnarray}
where $\chi_{1 - \alpha}^2(m)$ denotes the $1 - \alpha$ quantile of the
chi-squared distribution with $m$ degrees of freedom. Note that in
general $\mathcal{M}_1(d_n) \nsubseteq\mathcal{M}_2(d_n)$ and vice versa.
The confidence region $\mathcal{M}_2(d_n)$ can be viewed as a~global
measure, where the impact of single elements $\{|\widehat
{\theta}_{i} - \theta_{i}|\}_{1 \leq i \leq d_n}$ is
negligible, which in turn leads to suboptimal confidence regions for
single elements. In contrast, $\mathcal{M}_1(d_n)$ can be viewed as a~local measure where single elements have a~large impact, which clearly
leads to significantly tighter bounds for the single elements $\{
|\widehat{\theta}_{i} - \theta_{i}|\}_{1 \leq i
\leq
d_n}$. This is a~very important issue for so-called \textit{subset
autoregressive} models; see Remark~\ref{remautosubset}.

Theorem~\ref{thmyeconvrgextreme} not only can be used to construct
simultaneous confidence bands for the Yule--Walker estimators $\widehat
{\bolds\Theta}_{{d_n}}$, but also provides a~test for the degree of an
$\operatorname{AR}(q)$-process. To be more precise, for an $\operatorname{AR}(q)$-process $\{
X_k\}_{k \in\Z}$ satisfying the assumptions of Theorem~\ref{thmyeconvrgextreme}, we formulate the null hypothesis $\mathcal
{H}_0\dvtx q \leq q_0$, and the alternative $\mathcal{H}_A\dvtx q >
q_0$. Since
for any fixed $k \geq1$
\[
P\Bigl(a_n^{-1}\Bigl({\sqrt{n}\max_{1 \leq i \leq k}} |
(\widehat{\gamma}_{i,i}^* \widehat{\sigma}^2(d_n))^{-1/2}\widehat
{\theta
}_{i}| - b_{n}\Bigr) \leq z \Bigr) \rightarrow1
\]
as $n$ increases, it follows immediately from Theorem~\ref{thmyeconvrgextreme} that under $\mathcal{H}_0$ we have
\[
P\Bigl(a_n^{-1}\Bigl({\sqrt{n}\max_{q_0 + k \leq i \leq d_n}} |
(\widehat{\gamma}_{i,i}^* \widehat{\sigma}^2(d_n))^{-1/2}\widehat
{\theta
}_{i}| - b_{n}\Bigr) \leq z \Bigr) \rightarrow\exp(- e^{-z})
\]
for any fixed integer $k \geq1$, since we are assuming that $\theta_i
= 0$ for $i > q_0$.
Conversely, it is not hard to verify (see the proof of Theorem~\ref{thmestdistrib} for details) that the quantity
\[
a_n^{-1}\Bigl({\sqrt{n}\max_{q_0 + 1\leq i \leq d_n}} |
(\widehat{\gamma}_{i,i}^* \widehat{\sigma}^2(d_n))^{-1/2}\widehat
{\theta
}_{i}| - b_{n}\Bigr)
\]
explodes under the alternative $\mathcal{H}_A\dvtx q > q_0$. This can be
used to establish a~lower bound for the order $q$ or to test if the
order was chosen sufficiently large. This is particularly useful if $q$
is large compared to the sample size and the magnitude of~${\bolds\Theta
}_q$, in which case the AIC and related criteria sometimes heavily fail
to get near the true order. More details on this subject and examples
are given in Section~\ref{secnumerical}. Generally speaking, such
situations are often encountered in \textit{subset autoregressive}
models; see Remark~\ref{remautosubset}.

The above conclusions lead to the following family of estimators
$\widehat{q}_{z_n}^{(1)}$ for~$q$. Let $z_n$ be a~monotone sequence
that tends to infinity as $n$ increases. Then we define the estimator
%
%
\begin{equation}\label{defnest}
\widehat{q}_{z_n}^{(1)} = \min\Bigl\{q \in\N\big|
a_n^{-1}
\Bigl({\sqrt{n}\max_{q + 1 \leq i \leq d_n}} |(\widehat{\gamma}_{i,i}^*
\widehat{\sigma}^2(d_n))^{-1/2}\widehat{\theta}_{i}| -
b_{n}
\Bigr) \leq z_n \Bigr\}.\hspace*{-22pt}
\end{equation}

Using the above ideas, it is not hard to show that the estimators
$\widehat{q}_{z_n}^{(1)}$ are consistent if $z_n$ does not grow too
fast. In fact, under some more conditions imposed on the sequence
$z_n$, we can even derive the asymptotic distribution of the estimators.
%
%
\begin{ass}\label{assmainstrong}
In addition to Assumption~\ref{assmain}, suppose that:
\begin{itemize}
\item${\sum_{i = 1}^{\infty}} |\theta_i| < \infty$,
$
|\theta_n| = \OO((\log n)^{-2 - \eta})$, $\eta> 0$,
\item$\E(\exp(\lambda|\varepsilon_k|)) < \infty$, for some
$\lambda> 0$ and all $k \in\Z$,
\item$|\alpha_i| = \OO(i^{-\beta})$, $\beta
> 3/2$.
\end{itemize}
\end{ass}
%
%
\begin{theorem}\label{thmestdistrib}
Let $\{X_k\}_{k \in\Z}$ be an $\operatorname{AR}(q)$-process such that
Assumption~\ref{assmainstrong} is valid. Assume in addition that
$\inf
_{h}|\gamma_{h,h}^* | > 0$ and $z_n = \OO( \log n
)$. Then if $z_n \to\infty$, the estimator $\widehat{q}_{z_n}^{(1)}$
in~(\ref{defnest}) is consistent. Moreover, the following expansion
is valid:
\[
P\bigl(\widehat{q}_{z_n}^{(1)} = k + q \bigr) = \frac
{e^{-z_n}}{d_n} +
\OOO\biggl(\frac{e^{-z_n}}{d_n} +
d_n^{-z_n^2 + 1}\biggr)
\]
for $k \in\N$, $k = \OO(n^{\delta})$, $\delta< 1/7$.
\end{theorem}
%
%
\begin{rem}
The stronger conditions of Assumption~\ref{assmainstrong} are
necessary to control the rate of convergence in Theorem~\ref{thmyeconvrgextreme}, which in turn allows for the explicit
expansion given above. This, however, also leads to the more
restrictive bound $q + k = \OO(d_n) = \OO(n^{\delta
}
)$, $\delta< 1/7$; see also Remark~\ref{remincreasedelta}. If we are
only interested in establishing consistency, then we may drop these
more restrictive assumptions; see in particular Theorem~\ref{thmexpandK} below.
\end{rem}
%
%
\begin{rem}\label{remestcomp}
Theorem~\ref{thmestdistrib} yields that in some sense the
estimators~$\widehat{q}_{z_n}^{(1)}$ possess a~discrete uniform asymptotic
distribution, which leads to the surprising conclusion
\[
P\bigl(\widehat{q}_{z_n}^{(1)} = 1 + q \bigr) \approx P
\bigl(\widehat
{q}_{z_n}^{(1)} = 1000 + q \bigr).
\]
This fact\vspace*{2pt} can be explained by the maximum function in the definition
of~$\widehat{q}_{z_n}^{(1)}$, more precisely, due to the weak dependence
of the\vadjust{\goodbreak} Yule--Walker estimators~$\widehat{\bolds\Theta}_{{d_n}}$. The
maximum function essentially does not care at which index~$i$ the
boundary $z_n$ is exceeded, and this results in the uniform
distribution. It turns out (see Section~\ref{secnumerical}) that a~modified
version of the estimator~$\widehat{q}_{z_n}^{(1)}$ is a~very
efficient preliminary estimator that establishes a~decent lower bound.
\end{rem}

An asymptotic uniform-type distribution clearly is not a~desirable
property for an estimator. However, similarly to Akaike's method, we
can introduce a~penalty function and construct different yet also
consistent estimators for the order $q$.
To this end, for $x \in\R$ put $(x)^+ = \max(0,x)$ and let $\Upsilon
_{n,i} = a_n^{-1}(\sqrt{n}|(\widehat{\gamma}_{i,i}^*\times
\widehat
{\sigma}^2(d_n))^{-1/2}\widehat{\theta}_{i}| - b_{n})$. Then
we introduce a~new estimator $\widehat{q}_{z_n}^{(2)}$ as
\[
\widehat{q}_{z_n}^{(2)} = \argmin_{q \in\N}\Bigl\{\max_{q + 1
\leq i \leq d_n} \{(\Upsilon_{n,i} - z_n)^+ \} +
\log(1 + q)\Bigr\}.
\]
More generally, let $\F= (f_d)_{d \in\N}$ be a~collection
of continuous functions such that:
\begin{itemize}
\item$f_d$ is a~map from $\R^{d+2}$ to $\R$,
\item$f_d(0,\ldots,0,q,d) < f_d(0,\ldots,0,q + 1,d)$ for all $d,q \in\N$,
\item if $a_n, d_n \to\infty$ as $n$ increases, then
$f_{d_n}(\ldots,a_n,\ldots,q,d_n) \to\infty$ as $n$ increases,
regardless of
the values of the other coordinates.
\end{itemize}
Define
%
%
\begin{equation}
\widehat{q}_{z_n}^{(f)} = \argmin_{q \in\N} f_{d_n}
\bigl(0,\ldots,0,(\Upsilon_{n,q + 1} - z_n)^+,\ldots,(\Upsilon_{n,d_n} -
z_n)^+,q,d_n \bigr).\hspace*{-30pt}
\end{equation}
Then arguing as in the proof of Theorem~\ref{thmestdistrib} it can be
shown that this constitutes a~consistent estimator for the true value
$q$. For example, the following estimator
\[
\widehat{q}_{z_n}^{(3)} = \argmin_{q \in\N}\biggl\{\sum_{q + 1
\leq i \leq d_n} (\Upsilon_{n,i} - z_n)^+ + q\biggr\}
\]
satisfies the conditions above and is consistent.
%
%
\begin{rem}\label{remautosubset}
Note that instead of defining a~specific order $q$, one can also
consider a~special lag configuration, for example, ${\bolds\Theta}_q =
(\theta_1, \theta_2,0,\ldots,0,\allowbreak\theta_{10}, \theta
_{11},\ldots,\theta_q
)^T$. Such configurations are commonly referred to as \textit
{subset autoregressive models}; see, for instance,
\cite{brockwell2005,Mcclave,mcleod2008,tong1977,whittle1963} and the
references therein. The $\mathrm{AIC}(m)$ and especially related
consistent criteria have problems dealing with such \textit{subset
autoregressive models}, which can be seen as follows. By Hannan
\cite{hannantimeseries}, Chapter VI, we have for $m \in\N$
%
%
\begin{eqnarray}
\widetilde{\mathrm{AIC}}(m)n^{-1} &=& \log(\widehat{\sigma
}^2(m)
) + 2 n^{-1} C_n m \nonumber\\[-8pt]\\[-8pt]
&=& \log\widehat{\phi}_{n,0} + \sum_{i = 1}^m \log\bigl(1
- \widehat{\theta}_i^2(m)\bigr) + 2 n^{-1} C_n m.\nonumber
\end{eqnarray}
This shows that in case of \textit{subset autoregressive models}, the
penalty function $2 n^{-1} C_n m$ is too severe and should be replaced,
at least in theory, by $2 n^{-1} C_n\times \sum_{i = 1}^m \ind_{\{\theta_i
\neq0\}} $, since this is impossible in practice. Of course the same
problem arises if some of the $\{\theta_i\}_{1 \leq i \leq
q}$ are close to zero. A~maximum based estimator like $\widehat
{q}_{z_n}^{(1)}$ gets less effected, which is empirically confirmed in
Section~\ref{secnumerical}.
\end{rem}

An often encountered theoretical assumption for estimators related to
$\mathrm{AIC}(m)$ is that the parameter space for $q$ is finite; that
is, it is usually assumed in advance that $q \in\{0,\ldots,K\}$, where $K$
is ``chosen sufficiently large,'' but finite. In~\cite{hannanarp},
$K$ is allowed to increase with the sample size with unknown rate,
which was specified later by An et al.~\cite{hannan1982}. Note,
however, that for the estimators defined above we allow $K = K_n =
d_n$. Before extending this result, we give precise definitions of BIC,
HQC, MIC ($=$miscellaneous information criterion) and SIC, as the
literature does not seem to be very clear on this subject, in
particular in the case of the BIC and SIC. In the sequel, the following
definitions are used:
%
%
\begin{eqnarray}\label{defnbicusw}
\mathrm{BIC}(m) &=& \mathrm{SIC}(m) = \log\widehat{\sigma}^2(m) + m
n^{-1}\log n , \nonumber\\
\mathrm{MIC}(m) &=& \log\widehat{\sigma}^2(m) + m /2 n^{-1}\log n ,
\\
\mathrm{HQC}(m) &=& \log\widehat{\sigma}^2(m) + n^{-1} 2c m \log
\log
n,\qquad c > 1.\nonumber
\end{eqnarray}
This means that we use the same definitions for BIC and SIC
(asymptotically), which is the case mostly encountered in the
literature. The MIC differs from the BIC by the choice of the constant
$1/2$ that naturally leads to a~less parsimonious criterion, which
performs quite well in the examples given in Section~\ref{secnumerical}. Using some of the results of Section~\ref{secyulewalkerdetails} and~\ref{secgaussianapprox}, one may prove
the following.
%
%
\begin{theorem}\label{thmexpandK}
Assume that the conditions of Theorem~\ref{thmyeconvrgextreme} hold,
and additionally assume that $\inf_h|1 - \theta_h^2| > 0$.
Let $C_n$ be a~positive sequence such that:
\begin{itemize}
\item$\lim_n C_n (2 \log\log n)^{-1} > 1$, $C_n =
\OOO(n)$,
\item$ \log d_n \leq \OOO(C_n)$.
\end{itemize}
Then the estimators for the order $q$ defined as
\[
\widehat{q}_n^* = \argmin_{0 \leq m \leq d_n}\bigl(\log
\widehat{\sigma}^2(m) + n^{-1} C_n m \bigr)
\]
are consistent.
\end{theorem}
%
%
\begin{rem}
Note that condition ${\inf_h}|1 - \theta_h^2| > 0$ essentially
is already provided by the causality condition in Assumption~\ref{assmain}.
\end{rem}

Theorem~\ref{thmexpandK} thus implies the bounds $d_n \in\OO
\{\OOO(n^{1/2}), \OOO(n^{1/2}), \OOO(\log n)\}$ for BIC, MIC and HQC, and thus
significantly improves the bounds provided by An et al.
\cite{hannan1982} [BIC: $\OO(\log n)$, HQC: $\OO(\log\log n) $]. On the
other hand, the setting in An et al.~\cite{hannan1982} is more
general, and it is also shown that the estimators are strongly
consistent.

\section{Simulation and numerical results}\label{secnumerical}

In this section we will perform a~small simulation study to compare
some of the previously mentioned estimators. We will look at the
performance in case of $\operatorname{AR}(6)$, $\operatorname{AR}(12)$ and $\operatorname{AR}(24)$ processes. The sample
size $n$ satisfies $n \in\{125,250,500,1000\}$; as for the dimension
$d_n$, we chose the functions $d_n \in\{ 2 \log n, 4 \log n, 6 \log n\}
$, and rounded up the values. This implies that the parameter space $q
\in\{0,\ldots,K\}$ satisfies $K \in\{10,12,13,14\}$, $K \in\{
20,23,25,29\}$, $K \in\{29,34,38,42\}$. For reference, note that $\{
\lceil\sqrt{125}\rceil,\lceil\sqrt{250}\rceil,\lceil\sqrt {500}\rceil
,\lceil\sqrt{1000}\rceil\} = \{12,16,23,32\}$. To introduce the
estimators $\widehat{q}_{z_n}^{(4)}(d_n), \widehat
{q}_{z_n}^{(5)}(d_n)$, we require\vspace*{1pt} some additional
notation. For $1 \leq k \leq d_n$, define $\{\widehat{\gamma
}_{i,i}^*(k) \}_{1 \leq i \leq k}$ and $\{\widehat{\theta}_{i}(k)\}_{1
\leq i \leq k}$ via the usual relation
%
%
\begin{equation}
\widehat{\Theta}_k = \widehat{\bolds\Gamma}_{k}^{-1} \widehat{\Phi}_k.
\end{equation}
The estimators are now defined as
\begin{eqnarray*}
\widehat{q}_{z_n}^{(4)}(k) &=& \min\Bigl\{q \in\N\big|
a_n^{-1}\Bigl({\sqrt{n}\max_{q + 1 \leq i \leq k}} |(\widehat
{\gamma
}_{i,i}^*(k) \widehat{\sigma}^2(k))^{-1/2}\widehat{\theta
}_{i}(k)|
- b_{n}\Bigr) \leq z_n \Bigr\}, \\
\widehat{q}_{z_n}^{(5)}(d_n) &=&
\max
_{1 \leq k \leq d_n}\widehat{q}_{z_n}^{(4)}(k).
\end{eqnarray*}
Note that the definition of $a_n, b_n$ remains unchanged. This
modification significantly improves the performance in practice, which
is due to the following reason: if one just considers the estimator
$\widehat{q}_{z_n}^{(4)}(d_n)$ and hence \textit{only the equation}
$\widehat{\Theta}_{d_n} = \widehat{\bolds\Gamma}{}^{-1}_{d_n} \widehat
{\Phi
}_{d_n}$, the bias may be quite large since the estimate $\widehat
{\bolds
\Gamma}{}^{-1}_{d_n}$ is rather poor for larger $d_n$. Note that this is
also true when computing the AIC or related criteria, which is a~well-established fact in the literature (cf.
\cite{akaike1977,timeseriesbrockwell,hannanarp,hannantimeseries}).
Hence one may
expect that the
``maximum'' version~$\widehat{q}_{z_n}^{(5)}(d_n)$ outperforms its
counterpart $\widehat{q}_{z_n}^{(4)}(d_n)$, which is indeed the case in
the examples given below. The values for $z_n$ were chosen as $z_n \in
\{x_n,y_n\}$, where $x_n$ satisfies $a_n x_n + b_n = 2.71$ for $n \in\{
125,250\}$, $a_n x_n + b_n = 2.91$ for $n \in\{500,1000\}$. Similarly,
we have $a_n y_n + b_n = 3$ for $n \in\{125,250\}$, $a_n y_n + b_n =
3.2$ for $n \in\{500,1000\}$. This means that the estimators get less
parsimonious when $d_n$ increases. Of course an adaption to maintain
the same confidence level is possible, but the general picture remains
the same.

For the criteria AIC, BIC, HQC and MIC we use the definitions given in
(\ref{defnaik}) and~(\ref{defnbicusw}); in case of HQC we choose $c
= 1$, since, as pointed out by Hannan and
Quinn~\cite{hannanarp},\vadjust{\goodbreak}
``it would seem pedantic to choose values as $c = 1.01$.'' The following
modifications are also considered:
%
%
\begin{eqnarray}
\mathrm{AIC}(m)^* &=& \max\bigl\{\mathrm{AIC}(m), \widehat
{q}_{y_n}^{(5)}(d_n)\bigr\},\nonumber\\
\mathrm{BIC}(m)^* &=& \max\bigl\{\mathrm
{BIC}(m), \widehat{q}_{y_n}^{(5)}(d_n)\bigr\}, \nonumber\\[-8pt]\\[-8pt]
\mathrm{HQC}(m)^* &=& \max\bigl\{\mathrm{HQC}(m), \widehat
{q}_{y_n}^{(5)}(d_n)\bigr\},\nonumber\\
\mathrm{MIC}(m)^* &=& \max\bigl\{\mathrm
{MIC}(m), \widehat{q}_{y_n}^{(5)}(d_n)\bigr\}.\nonumber
\end{eqnarray}

All simulations were carried out using the
program $R$;\footnote{%
\texttt{\href{http://portal.tugraz.at/portal/page/portal/TU\_Graz/Einrichtungen/Institute/Homepages/i5060/research/R\_Code}{http://portal.tugraz.at/portal/page/portal/TU\_Graz/Einrichtungen/Institute/}
\href{http://portal.tugraz.at/portal/page/portal/TU\_Graz/Einrichtungen/Institute/Homepages/i5060/research/R\_Code}{Homepages/i5060/research/R\_Code}}.}
in order to get a~sample of size $n$, a~sample path of size $1000 + n$
was produced and the first $1000$ observations were discarded.

Generally speaking, unreported simulations show that in many cases the
modified criteria $\mathrm{AIC}(m)^*,\mathrm{BIC}(m)^*,\ldots$
perform
nearly identically as the nonmodified ones $\mathrm{AIC}(m),\mathrm
{BIC}(m),\ldots.$ This is in particular the case when dealing with full
parameter sets, that is, $\theta_i \neq0$, $1 \leq i \leq q$, and
$\theta_q$ is sufficiently large. If this is the case, the performance
of the estimators $\widehat{q}_{x_n}^{(5)}(d_n)$, $\widehat{q}_{y_n}^{(5)}(d_n)$
is somewhere between the $\mathrm{BIC}(m)$ and
$\mathrm{HQC}(m)$. On the other hand, if the model is not full and/or
the order $q$ is sufficiently large, then the differences can be quite
striking. The aim of the following examples is to illustrate this behavior.

\subsection{$\operatorname{AR}(6)$}

First note that the definitions of $x_n, y_n$ result in
\begin{eqnarray*}
P({\max}|\xiv| \leq2.71 ) &\geq& 0.92,\qquad P(\max
|\xiv|
\leq3 ) \geq0.97,\qquad d_n \in\{10,12\},\\
P({\max}|\xiv| \leq2.91 ) &\geq& 0.95,\qquad P(\max
|\xiv|
\leq3.2 ) \geq0.98,\qquad d_n \in\{13,14\},
\end{eqnarray*}
where $\xiv= (\xi_1,\ldots,\xi_{d_n})^T$ is a~$d_n$-dimensional
mean-zero Gaussian random vector where the covariance matrix is the identity.

The results shown in Tables~\ref{tabquartielsaro} and~\ref{tabquartielsaro1} hint at what is to be expected in case of full
models, namely that the modifications $\mathrm{AIC}(m)^*,\mathrm
{BIC}(m)^*,\ldots$ perform nearly as well as the normal versions
$\mathrm
{AIC}(m),\mathrm{BIC}(m),\ldots.$ The estimators $\widehat
{q}_{x_n}^{(5)}(d_n)$, $\widehat{q}_{y_n}^{(5)}(d_n)$ perform also quite well.

%
%
\begin{table}
\def\arraystretch{0.95}
\caption{Simulation of an $\operatorname{AR}(6)$ process with coefficients $\Theta_6 =
(0.1,-0.3,0.05,0.2,-0.1,0.2)^T$, $\varepsilon\sim\Gaussian(0,1)$, 1000
repetitions, $d_n \in\{10,12\}$}
\label{tabquartielsaro}
\begin{tabular*}{\tablewidth}{@{\extracolsep{\fill}}lrrrrrrrrrrr@{}}
\hline
$\bolds{n}$& \multicolumn{1}{c}{$\widehat{\bolds{q}}$} &
\multicolumn{1}{c}{\textbf{AIC}\hspace*{-1pt}}
& \multicolumn{1}{c}{\textbf{AIC*}\hspace*{-3pt}}
& \multicolumn{1}{c}{\textbf{BIC}\hspace*{-1pt}}
& \multicolumn{1}{c}{\textbf{BIC*}\hspace*{-3pt}}
& \multicolumn{1}{c}{\textbf{HQC}\hspace*{-3pt}}
& \multicolumn{1}{c}{\textbf{HQC*}\hspace*{-5pt}}
& \multicolumn{1}{c}{\textbf{MIC}\hspace*{-2pt}}
& \multicolumn{1}{c}{\textbf{MIC*}\hspace*{-4pt}}
& \multicolumn{1}{c}{$\bolds{\widehat{q}_{y_n}^{(5)}}$\hspace*{-2pt}}
& \multicolumn{1}{c@{}}{$\bolds{\widehat{q}_{x_n}^{(5)}}$} \\
\hline
125 & $< $5 & 428 & 427 & 943 & 808 & 746 & 704 & 550 & 545 & 816 &
701\hspace*{2pt}\\
& 5 & 65 & 65 & 10 & 30 & 32 & 40 & 58 & 58 & 28 & 41\hspace*{2pt}\\
& \textbf{6} & \textbf{344} & \textbf{341} & \textbf{45} & \textbf{143} & \textbf{191} & \textbf{214} &
\textbf{295}
& \textbf{294} & \textbf{137} & \textbf{196}\hspace*{2pt}\\
& 7 & 66 & 65 & 1 & 5 & 23 & 24 & 54 & 53 & 5 & 14\hspace*{2pt}\\
& $< $7 & 97 & 102 & 1 & 14 & 8 & 18 & 43 & 50 & 14 & 48\hspace*{2pt}\\
[6pt]
250 & $< $5 & 93 & 89 & 693 & 432 & 328 & 282 & 202 & 188 & 440 &
299\hspace*{2pt}\\
& 5 & 24 & 23 & 14 & 32 & 32 & 32 & 33 & 31 & 42 & 38\hspace*{2pt}\\
& \textbf{6} & \textbf{646} & \textbf{632} & \textbf{287} & \textbf{481} & \textbf{586} & \textbf{595} &
\textbf{649}
& \textbf{634} & \textbf{467} & \textbf{543}\hspace*{2pt}\\
& 7 & 96 & 95 & 5 & 8 & 37 & 35 & 74 & 73 & 4 & 9\hspace*{2pt}\\
& $> $7 & 141 & 161 & 1 & 47 & 17 & 56 & 42 & 74 & 47 & 111\hspace*{2pt}\\
\hline
\end{tabular*}
\vspace*{-4pt}
\end{table}

%
%
\begin{table}[b]
\def\arraystretch{0.95}
\vspace*{-4pt}
\caption{Simulation of an $\operatorname{AR}(6)$ process with coefficients $\Theta_6 =
(0.1,-0.3,0.05,0.2,-0.1,0.2)^T$, $\varepsilon\sim\Gaussian(0,1)$, 1000
repetitions, $d_n \in\{13,14\}$}
\label{tabquartielsaro1}
\begin{tabular*}{\tablewidth}{@{\extracolsep{\fill}}lrrrrrrrrrrr@{}}
\hline
$\bolds{n}$& \multicolumn{1}{c}{$\widehat{\bolds{q}}$} &
\multicolumn{1}{c}{\textbf{AIC}\hspace*{-1pt}}
& \multicolumn{1}{c}{\textbf{AIC*}\hspace*{-3pt}}
& \multicolumn{1}{c}{\textbf{BIC}\hspace*{-1pt}}
& \multicolumn{1}{c}{\textbf{BIC*}\hspace*{-3pt}}
& \multicolumn{1}{c}{\textbf{HQC}\hspace*{-3pt}}
& \multicolumn{1}{c}{\textbf{HQC*}\hspace*{-5pt}}
& \multicolumn{1}{c}{\textbf{MIC}\hspace*{-2pt}}
& \multicolumn{1}{c}{\textbf{MIC*}\hspace*{-4pt}}
& \multicolumn{1}{c}{$\bolds{\widehat{q}_{y_n}^{(5)}}$\hspace*{-2pt}}
& \multicolumn{1}{c@{}}{$\bolds{\widehat{q}_{x_n}^{(5)}}$} \\
\hline
\hphantom{0}500 & $< $5 & 1 & 1 & 177 & 75 & 29 & 25 & 15 & 15 & 86 & 52\hspace*{2pt}\\
& 5 & 3 & 3 & 9 & 11 & 6 & 6 & 3 & 3 & 17 & 14\hspace*{2pt}\\
& \textbf{6} & \textbf{730} & \textbf{713} & \textbf{805} & \textbf{874} & \textbf{913} & \textbf{889} &
\textbf{892}
& \textbf{867} & \textbf{865} & \textbf{849}\hspace*{2pt}\\
& 7 & 108 & 108 & 8 & 8 & 42 & 42 & 57 & 57 & 0 & 2\hspace*{2pt}\\
& $< $7 & 158 & 175 & 1 & 32 & 10 & 38 & 33 & 58 & 32 & 83\hspace*{2pt}\\
[6pt]
1000 & $< $5 & 0 & 0 & 3 & 0 & 0 & 0 & 0 & 0 & 0 & 0\hspace*{2pt}\\
& 5 & 0 & 0 & 0 & 0 & 0 & 0 & 0 & 0 & 0 & 0\hspace*{2pt}\\
& \textbf{6} & \textbf{724} & \textbf{709} & \textbf{990} & \textbf{951} & \textbf{952} & \textbf{917} &
\textbf{934}
& \textbf{901} & \textbf{955} & \textbf{885}\hspace*{2pt}\\
& 7 & 103 & 101 & 7 & 9 & 36 & 34 & 47 & 44 & 5 & 7\hspace*{2pt}\\
& $> $7 & 173 & 190 & 0 & 40 & 12 & 49 & 19 & 55 & 40 & 108\hspace*{2pt}\\
\hline
\end{tabular*}
\end{table}

%
%
\begin{table}
\caption{Simulation of an $\operatorname{AR}(6)$ process with coefficients $\Theta_6 =
(0.1,0,0.05,0,0,0.2)^T$, $\varepsilon\sim\Gaussian(0,1)$, 1000~repetitions, $d_n \in\{10,12\}$}
\label{tabquartiels1}
\begin{tabular*}{\tablewidth}{@{\extracolsep{\fill}}lrrrrrrrrrrr@{}}
\hline
$\bolds{n}$& \multicolumn{1}{c}{$\widehat{\bolds{q}}$} &
\multicolumn{1}{c}{\textbf{AIC}\hspace*{-1pt}}
& \multicolumn{1}{c}{\textbf{AIC*}\hspace*{-3pt}}
& \multicolumn{1}{c}{\textbf{BIC}\hspace*{-1pt}}
& \multicolumn{1}{c}{\textbf{BIC*}\hspace*{-3pt}}
& \multicolumn{1}{c}{\textbf{HQC}\hspace*{-3pt}}
& \multicolumn{1}{c}{\textbf{HQC*}\hspace*{-5pt}}
& \multicolumn{1}{c}{\textbf{MIC}\hspace*{-2pt}}
& \multicolumn{1}{c}{\textbf{MIC*}\hspace*{-4pt}}
& \multicolumn{1}{c}{$\bolds{\widehat{q}_{y_n}^{(5)}}$\hspace*{-2pt}}
& \multicolumn{1}{c@{}}{$\bolds{\widehat{q}_{x_n}^{(5)}}$} \\
\hline
125 & $<$5 & 719 & 699 & 998 & 854 & 944 & 842 & 839 & 787 & 854 &
747\hspace*{2pt}\\
& 5& 11 & 11 & 0 & 0 & 2 & 2 & 7 & 7 & 0 & 11\hspace*{2pt}\\
& \textbf{6} & \textbf{168} & \textbf{181} & \textbf{2} & \textbf{124} & \textbf{43} & \textbf{126} &
\textbf{107} &
\textbf{145} & \textbf{124} & \textbf{184}\hspace*{2pt}\\
& 7 & 44 & 44 & 0 & 4 & 8 & 11 & 23 & 24 & 4 & 8\hspace*{2pt}\\
& $< $7 & 58 & 65 & 0 & 18 & 3 & 19 & 24 & 37 & 18 & 50 \hspace*{2pt}\\
[6pt]
250 & $<$5 & 290 & 276 & 960 & 437 & 723 & 424 & 550 & 396 & 438 &
321\hspace*{2pt}\\
& 5 & 6 & 6 & 0 & 3 & 2 & 3 & 5 & 5 & 3 & 5\hspace*{2pt}\\
& \textbf{6} & \textbf{491} & \textbf{488} & \textbf{39} & \textbf{513} & \textbf{245} & \textbf{503} &
\textbf{376} & \textbf{494} & \textbf{513} & \textbf{573}\hspace*{2pt}\\
& 7 & 91 & 90 & 1 & 2 & 21 & 21 & 40 & 40 & 1 & 7\hspace*{2pt}\\
& $>$7 & 122 & 140 & 0 & 45 & 9 & 49 & 29 & 65 & 45 & 94\hspace*{2pt}\\
\hline
\end{tabular*}
\end{table}

%
%
\begin{table}[b]
\caption{Simulation of an $\operatorname{AR}(6)$ process with coefficients $\Theta_6 =
(0.1,0,0.05,0,0,0.2)^T$, $\varepsilon\sim\Gaussian(0,1)$, 1000~repetitions, $d_n \in\{13,14\}$}
\label{tabquartiels11}
\begin{tabular*}{\tablewidth}{@{\extracolsep{\fill}}lrrrrrrrrrrr@{}}
\hline
$\bolds{n}$& \multicolumn{1}{c}{$\widehat{\bolds{q}}$} &
\multicolumn{1}{c}{\textbf{AIC}\hspace*{-1pt}}
& \multicolumn{1}{c}{\textbf{AIC*}\hspace*{-3pt}}
& \multicolumn{1}{c}{\textbf{BIC}\hspace*{-1pt}}
& \multicolumn{1}{c}{\textbf{BIC*}\hspace*{-3pt}}
& \multicolumn{1}{c}{\textbf{HQC}\hspace*{-3pt}}
& \multicolumn{1}{c}{\textbf{HQC*}\hspace*{-5pt}}
& \multicolumn{1}{c}{\textbf{MIC}\hspace*{-2pt}}
& \multicolumn{1}{c}{\textbf{MIC*}\hspace*{-4pt}}
& \multicolumn{1}{c}{$\bolds{\widehat{q}_{y_n}^{(5)}}$\hspace*{-2pt}}
& \multicolumn{1}{c@{}}{$\bolds{\widehat{q}_{x_n}^{(5)}}$} \\
\hline
\hphantom{0}500 & $< $5 & 21 & 21 & 761 & 102 & 267 & 98 & 164 & 85 & 102 & 56\hspace*{2pt}\\
& 5 & 0 & 0 & 1 & 0 & 0 & 0 & 0 & 0 & 0 & 1\hspace*{2pt}\\
& \textbf{6} & \textbf{663} & \textbf{655} & \textbf{234} & \textbf{871} & \textbf{675} & \textbf{822} &
\textbf{736}
& \textbf{796} & \textbf{874} & \textbf{863}\hspace*{2pt}\\
& 7 & 125 & 124 & 4 & 3 & 50 & 49 & 69 & 68 & 0 & 10\hspace*{2pt}\\
& $< $7 & 191 & 200 & 0 & 24 & 8 & 31 & 31 & 51 & 24 & 70\hspace*{2pt}\\
[6pt]
1000 & $< $5 & 0 & 0 & 168 & 1 & 3 & 1 & 1 & 1 & 1 & 0\hspace*{2pt}\\
& 5 & 0 & 0 & 0 & 0 & 0 & 0 & 0 & 0 & 0 & 0\hspace*{2pt}\\
& \textbf{6} & \textbf{702} & \textbf{683} & \textbf{822} & \textbf{949} & \textbf{940} & \textbf{905} &
\textbf{919}
& \textbf{887} & \textbf{955} & \textbf{898}\hspace*{2pt}\\
& 7 & 121 & 119 & 9 & 9 & 43 & 42 & 52 & 52 & 3 & 9\hspace*{2pt}\\
& $> $7 & 177 & 198 & 1 & 41 & 14 & 52 & 28 & 60 & 41 & 93\hspace*{2pt}\\
\hline
\end{tabular*}
\end{table}

Contrary to the previous results, Tables~\ref{tabquartiels1} and
\ref
{tabquartiels11} show the difference of the modified estimators [and
$\widehat{q}_{x_n}^{(5)}(d_n), \widehat{q}_{y_n}^{(5)}(d_n)$], if the
model is very sparse. Except for the case $n = 1000$, the modifications
are notably better.

\subsection{$\operatorname{AR}(12)$}

The definitions of $x_n, y_n$ result in
\begin{eqnarray*}
P({\max}|\xiv| \leq2.71 ) &\geq&0.85,\qquad
P({\max}|\xiv|
\leq3 ) \geq0.94, \qquad d_n \in\{20,23\},\\
P({\max}|\xiv| \leq2.91 ) &\geq&0.9,\qquad P({\max}
|\xiv|\leq3.2 ) \geq0.96, \qquad d_n \in\{25,29\},
\end{eqnarray*}
where $\xiv= (\xi_1,\ldots,\xi_{d_n})^T$ is a~$d_n$-dimensional
mean-zero Gaussian random vector where the covariance matrix is the identity.

%
%
\begin{table}
\caption{Simulation of an $\operatorname{AR}(12)$ process with nonzero coefficients
$\theta_1 = 0.1$, $\theta_3 = -0.4$, $\theta_5 = 0.5$, $\theta_7 =
-0.1$, $\theta_8 = 0.05$, $\theta_{10} = -0.3$, $\theta_{12} = 0.2$,
$\varepsilon\sim\Gaussian(0,1)$, 1000 repetitions, $d_n \in\{20,23\}$}
\label{tabquartiels2}
\begin{tabular*}{\tablewidth}{@{\extracolsep{\fill}}lrrrrrrrrrrr@{}}
\hline
$\bolds{n}$& \multicolumn{1}{c}{$\widehat{\bolds{q}}$} &
\multicolumn{1}{c}{\textbf{AIC}\hspace*{-1pt}}
& \multicolumn{1}{c}{\textbf{AIC*}\hspace*{-3pt}}
& \multicolumn{1}{c}{\textbf{BIC}\hspace*{-1pt}}
& \multicolumn{1}{c}{\textbf{BIC*}\hspace*{-3pt}}
& \multicolumn{1}{c}{\textbf{HQC}\hspace*{-3pt}}
& \multicolumn{1}{c}{\textbf{HQC*}\hspace*{-5pt}}
& \multicolumn{1}{c}{\textbf{MIC}\hspace*{-2pt}}
& \multicolumn{1}{c}{\textbf{MIC*}\hspace*{-4pt}}
& \multicolumn{1}{c}{$\bolds{\widehat{q}_{y_n}^{(5)}}$\hspace*{-2pt}}
& \multicolumn{1}{c@{}}{$\bolds{\widehat{q}_{x_n}^{(5)}}$} \\
\hline
125 & $<$11 & 705 & 701 & 995 & 966 & 931 & 917 & 812 & 807 & 969 &
929\hspace*{2pt}\\
& 11 & 79 & 79 & 2 & 3 & 22 & 22 & 54 & 54 & 1 & 2\hspace*{2pt}\\
& \textbf{12} & \textbf{141} & \textbf{141} & \textbf{3} & \textbf{23} & \textbf{40} & \textbf{47} &
\textbf{97} & \textbf{98} & \textbf{22} & \textbf{47}\hspace*{2pt}\\
& 13 & 48 & 48 & 0 & 4 & 6 & 9 & 30 & 30 & 4 & 11\hspace*{2pt}\\
& $>$13 & 27 & 31 & 0 & 4 & 1 & 5 & 7 & 11 & 4 & 11\hspace*{2pt}\\
[6pt]
250 & $<$11 & 257 & 257 & 854 & 730 & 573 & 560 & 423 & 421 & 748 &
620\hspace*{2pt}\\
& 11 & 39 & 39 & 9 & 10 & 31 & 31 & 39 & 39 & 3 & 11\hspace*{2pt}\\
& \textbf{12} & \textbf{495} & \textbf{493} & \textbf{135} & \textbf{247} & \textbf{349} & \textbf{356} &
\textbf{442}
& \textbf{441} & \textbf{237} & \textbf{313}\hspace*{2pt}\\
& 13 & 115 & 115 & 2 & 4 & 40 & 40 & 65 & 65 & 3 & 13\hspace*{2pt}\\
& $>$13 & 94 & 96 & 0 & 9 & 7 & 13 & 31 & 34 & 9 & 43\hspace*{2pt}\\
\hline
\end{tabular*}
\end{table}

%
%
\begin{table}[b]
\caption{Simulation of an $\operatorname{AR}(12)$ process with nonzero coefficients
$\theta_1 = 0.1$, $\theta_3 = -0.4$, $\theta_5 = 0.5$, $\theta_7 =
-0.1$, $\theta_8 = 0.05$, $\theta_{10} = -0.3$, $\theta_{12} = 0.2$,
$\varepsilon\sim\Gaussian(0,1)$, 1000 repetitions, $d_n \in\{25,28\}$}
\label{tabquartiels22}
\begin{tabular*}{\tablewidth}{@{\extracolsep{\fill}}lrrrrrrrrrrr@{}}
\hline
$\bolds{n}$& \multicolumn{1}{c}{$\widehat{\bolds{q}}$} &
\multicolumn{1}{c}{\textbf{AIC}\hspace*{-1pt}}
& \multicolumn{1}{c}{\textbf{AIC*}\hspace*{-3pt}}
& \multicolumn{1}{c}{\textbf{BIC}\hspace*{-1pt}}
& \multicolumn{1}{c}{\textbf{BIC*}\hspace*{-3pt}}
& \multicolumn{1}{c}{\textbf{HQC}\hspace*{-3pt}}
& \multicolumn{1}{c}{\textbf{HQC*}\hspace*{-5pt}}
& \multicolumn{1}{c}{\textbf{MIC}\hspace*{-2pt}}
& \multicolumn{1}{c}{\textbf{MIC*}\hspace*{-4pt}}
& \multicolumn{1}{c}{$\bolds{\widehat{q}_{y_n}^{(5)}}$\hspace*{-2pt}}
& \multicolumn{1}{c@{}}{$\bolds{\widehat{q}_{x_n}^{(5)}}$} \\
\hline
\hphantom{0}500 & $<$11 & 19 & 19 & 367 & 256 & 110 & 106 & 75 & 73 & 269 & 183\hspace*{2pt}\\
& 11 & 4 & 4 & 4 & 4 & 6 & 6 & 6 & 6 & 2 & 2\hspace*{2pt}\\
& \textbf{12} & \textbf{684} & \textbf{680} & \textbf{618} & \textbf{705} & \textbf{808} & \textbf{793} &
\textbf{808} & \textbf{797} & \textbf{702} & \textbf{758}\hspace*{2pt}\\
& 13 & 129 & 128 & 10 & 12 & 63 & 62 & 78 & 76 & 4 & 8\hspace*{2pt}\\
& $>$13 & 164 & 169 & 1 & 23 & 13 & 33 & 33 & 48 & 23 & 49\hspace*{2pt}\\
[6pt]
1000 & $<$11 & 0 & 0 & 11 & 2 & 0 & 0 & 0 & 0 & 2 & 1\hspace*{2pt}\\
& 11 & 0 & 0 & 0 & 0 & 0 & 0 & 0 & 0 & 0 & 0\hspace*{2pt}\\
& \textbf{12} & \textbf{679} & \textbf{676} & \textbf{970} & \textbf{947} & \textbf{925} & \textbf{900} &
\textbf{896} & \textbf{873} & \textbf{958} & \textbf{914}\hspace*{2pt}\\
& 13 & 151 & 150 & 17 & 17 & 61 & 60 & 79 & 78 & 6 & 13\hspace*{2pt}\\
& $>$13 & 170 & 174 & 2 & 34 & 14 & 40 & 25 & 49 & 34 & 72\hspace*{2pt}\\
\hline
\end{tabular*}
\end{table}

%
%
\begin{table}
\caption{Simulation of an $\operatorname{AR}(12)$ process with nonzero coefficients
$\theta_1 = 0.1$, $\theta_3 = -0.4$, $\theta_{12} = 0.2$, $\varepsilon
\sim\Gaussian(0,1)$, 1000 repetitions, $d_n \in\{20,23\}$}
\label{tabquartiels3}
\begin{tabular*}{\tablewidth}{@{\extracolsep{\fill}}lrrrrrrrrrrr@{}}
\hline
$\bolds{n}$& \multicolumn{1}{c}{$\widehat{\bolds{q}}$} &
\multicolumn{1}{c}{\textbf{AIC}\hspace*{-1pt}}
& \multicolumn{1}{c}{\textbf{AIC*}\hspace*{-3pt}}
& \multicolumn{1}{c}{\textbf{BIC}\hspace*{-1pt}}
& \multicolumn{1}{c}{\textbf{BIC*}\hspace*{-3pt}}
& \multicolumn{1}{c}{\textbf{HQC}\hspace*{-3pt}}
& \multicolumn{1}{c}{\textbf{HQC*}\hspace*{-5pt}}
& \multicolumn{1}{c}{\textbf{MIC}\hspace*{-2pt}}
& \multicolumn{1}{c}{\textbf{MIC*}\hspace*{-4pt}}
& \multicolumn{1}{c}{$\bolds{\widehat{q}_{y_n}^{(5)}}$\hspace*{-2pt}}
& \multicolumn{1}{c@{}}{$\bolds{\widehat{q}_{x_n}^{(5)}}$} \\
\hline
125 & $<$10 & 884 & 853 & 1000 & 920 & 995 & 920 & 963 & 910 & 920 &
861\hspace*{2pt}\\
& 11 & 3 & 3 & 0 & 0 & 0 & 0 & 1 & 1 & 0 & 3\hspace*{2pt}\\
& \textbf{12} & \textbf{68} & \textbf{94} & \textbf{0} & \textbf{71} & \textbf{5} & \textbf{71} & \textbf{25} &
\textbf{70}
& \textbf{71} & \textbf{114}\hspace*{2pt}\\
& 13 & 11 & 13 & 0 & 3 & 0 & 3 & 4 & 7 & 3 & 5\hspace*{2pt}\\
& $>$13 & 34 & 37 & 0 & 6 & 0 & 6 & 7 & 12 & 6 & 17\hspace*{2pt}\\
[6pt]
250 & $<$10 & 509 & 421 & 999 & 555 & 934 & 552 & 792 & 530 & 555 &
424\hspace*{2pt}\\
& 11 & 3 & 3 & 0 & 3 & 0 & 2 & 2 & 3 & 3 & 4\hspace*{2pt}\\
& \textbf{12} & \textbf{340} & \textbf{419} & \textbf{1} & \textbf{421} & \textbf{59} & \textbf{419} &
\textbf{170} &
\textbf{416} & \textbf{421} & \textbf{514}\hspace*{2pt}\\
& 13 & 67 & 68 & 0 & 2 & 4 & 6 & 18 & 19 & 2 & 5\hspace*{2pt}\\
& $>$13 & 81 & 89 & 0 & 19 & 3 & 21 & 18 & 32 & 19 & 53\hspace*{2pt}\\
\hline
\end{tabular*}\vspace*{6pt}
\end{table}

%
%
\begin{table}[b]
\vspace*{6pt}
\caption{Simulation of an $\operatorname{AR}(12)$ process with nonzero coefficients
$\theta_1 = 0.1$, $\theta_3 = -0.4$, $\theta_{12} = 0.2$, $\varepsilon
\sim\Gaussian(0,1)$, 1000 repetitions, $d_n \in\{25,28\}$}
\label{tabquartiels33}
\begin{tabular*}{\tablewidth}{@{\extracolsep{\fill}}lrrrrrrrrrrr@{}}
\hline
$\bolds{n}$& \multicolumn{1}{c}{$\widehat{\bolds{q}}$} &
\multicolumn{1}{c}{\textbf{AIC}\hspace*{-1pt}}
& \multicolumn{1}{c}{\textbf{AIC*}\hspace*{-3pt}}
& \multicolumn{1}{c}{\textbf{BIC}\hspace*{-1pt}}
& \multicolumn{1}{c}{\textbf{BIC*}\hspace*{-3pt}}
& \multicolumn{1}{c}{\textbf{HQC}\hspace*{-3pt}}
& \multicolumn{1}{c}{\textbf{HQC*}\hspace*{-5pt}}
& \multicolumn{1}{c}{\textbf{MIC}\hspace*{-2pt}}
& \multicolumn{1}{c}{\textbf{MIC*}\hspace*{-4pt}}
& \multicolumn{1}{c}{$\bolds{\widehat{q}_{y_n}^{(5)}}$\hspace*{-2pt}}
& \multicolumn{1}{c@{}}{$\bolds{\widehat{q}_{x_n}^{(5)}}$} \\
\hline
\hphantom{0}500 & $<$11 & 77 & 58 & 983 & 125 & 613 & 125 & 402 & 115 & 125 & 78\hspace*{2pt}\\
& 11 & 0 & 0 & 0 & 2 & 0 & 2 & 0 & 1 & 2 & 1\hspace*{2pt}\\
& \textbf{12} & \textbf{663} & \textbf{678} & \textbf{17} & \textbf{858} & \textbf{360} & \textbf{834} &
\textbf{532}
& \textbf{808} & \textbf{858} & \textbf{870}\hspace*{2pt}\\
& 13 & 104 & 103 & 0 & 3 & 15 & 16 & 39 & 40 & 3 & 4\hspace*{2pt}\\
& $>$13 & 156 & 161 & 0 & 12 & 12 & 23 & 27 & 36 & 12 & 47\hspace*{2pt}\\
[6pt]
1000 & $< $11 & 0 & 0 & 689 & 2 & 67 & 2 & 35 & 2 & 2 & 2\hspace*{2pt}\\
& 11 & 0 & 0 & 0 & 0 & 0 & 0 & 0 & 0 & 0 & 0\hspace*{2pt}\\
& \textbf{12} & \textbf{706} & \textbf{701} & \textbf{307} & \textbf{971} & \textbf{880} & \textbf{926} &
\textbf{893}
& \textbf{907} & \textbf{972} & \textbf{936}\hspace*{2pt}\\
& 13 & 124 & 123 & 2 & 2 & 39 & 38 & 54 & 53 & 1 & 3\hspace*{2pt}\\
& $>$13 & 170 & 176 & 2 & 25 & 14 & 34 & 18 & 38 & 25 & 59\hspace*{2pt}\\
\hline
\end{tabular*}
\end{table}

The results are depicted in Tables~\ref{tabquartiels2},~\ref{tabquartiels22},~\ref{tabquartiels3} and~\ref{tabquartiels33},
and are quite similar to the case of the $\operatorname{AR}(6)$ processes. If the model
is rather full, $\mathrm{AIC}(m)^*,\mathrm{BIC}(m)^*,\ldots$ perform nearly
as well as the normal versions $\mathrm{AIC}(m),\mathrm{BIC}(m),\ldots,$
whereas in case of the sparse model, a~significant difference can be observed.\vspace*{-2pt}

\subsection{$\operatorname{AR}(24)$}

In this case, the definitions of $x_n, y_n$ result in
\begin{eqnarray*}
P({\max}|\xiv| \leq2.71 ) &\geq&0.795,\qquad P({\max}
|\xiv|
\leq3 ) \geq0.912,\qquad d_n \in\{29,34\},\\
P({\max}|\xiv| \leq2.91 ) &\geq&0.86,\qquad P({\max}
|\xiv|
\leq3.2 ) \geq0.94,\qquad d_n \in\{38,42\},
\end{eqnarray*}
where $\xiv= (\xi_1,\ldots,\xi_{d_n})^T$ is a~$d_n$-dimensional
mean-zero Gaussian random vector where the covariance matrix\vadjust{\goodbreak} is the
%
%
\begin{table}
\def\arraystretch{0.9}
\caption{Simulation of an $\operatorname{AR}(24)$ process with nonzero coefficients
$\theta_1 = 0.6$, $\theta_2 = -0.1$, $\theta_4 = 0.05$, $\theta_7 =
0.15$, $\theta_8 = -0.27$, $\theta_{10} = 0.1$, $\theta_{12} = -0.2$,
$\theta_{15} = -0.25$, $\theta_{18} = 0.05$, $\theta_{20} = 0.1$,
$\theta_{21} = -0.3$, $\theta_{24} = 0.17$, $\varepsilon\sim\Gaussian
(0,1)$, 1000 repetitions, $d_n \in\{29,34\}$}
\label{tabquartiels4}
\begin{tabular*}{\tablewidth}{@{\extracolsep{\fill}}lrrrrrrrrrrr@{}}
\hline
$\bolds{n}$& \multicolumn{1}{c}{$\widehat{\bolds{q}}$} &
\multicolumn{1}{c}{\textbf{AIC}\hspace*{-1pt}}
& \multicolumn{1}{c}{\textbf{AIC*}\hspace*{-3pt}}
& \multicolumn{1}{c}{\textbf{BIC}\hspace*{-1pt}}
& \multicolumn{1}{c}{\textbf{BIC*}\hspace*{-3pt}}
& \multicolumn{1}{c}{\textbf{HQC}\hspace*{-3pt}}
& \multicolumn{1}{c}{\textbf{HQC*}\hspace*{-5pt}}
& \multicolumn{1}{c}{\textbf{MIC}\hspace*{-2pt}}
& \multicolumn{1}{c}{\textbf{MIC*}\hspace*{-4pt}}
& \multicolumn{1}{c}{$\bolds{\widehat{q}_{y_n}^{(5)}}$\hspace*{-2pt}}
& \multicolumn{1}{c@{}}{$\bolds{\widehat{q}_{x_n}^{(5)}}$} \\
\hline
125 & $<$23 & 972 & 970 & 1000 & 996 & 1000 & 996 & 992 & 990 & 996 &
989\hspace*{2pt}\\
& 23 & 12 & 12 & 0 & 1 & 0 & 1 & 5 & 5 & 1 & 2\hspace*{2pt}\\
& \textbf{24} & \textbf{3} & \textbf{3} & \textbf{0} & \textbf{1} & \textbf{0} & \textbf{1} & \textbf{1} &
\textbf{1} & \textbf{1} & \textbf{6}\hspace*{2pt}\\
& 25 & 10 & 10 & 0 & 0 & 0 & 0 & 2 & 2 & 0 & 1\hspace*{2pt}\\
& $>$25 & 3 & 5 & 0 & 2 & 0 & 2 & 0 & 2 & 2 & 2\hspace*{2pt}\\
[6pt]
250 & $<$23 & 518 & 516 & 995 & 923 & 872 & 840 & 727 & 717 & 924 &
845\hspace*{2pt}\\
& 23 & 120 & 120 & 2 & 13 & 48 & 50 & 77 & 78 & 12 & 25\hspace*{2pt}\\
& \textbf{24} & \textbf{185} & \textbf{186} & \textbf{3} & \textbf{57} & \textbf{67} & \textbf{90} &
\textbf{135} & \textbf{138} & \textbf{57} & \textbf{98}\hspace*{2pt}\\
& 25 & 89 & 89 & 0 & 1 & 7 & 8 & 38 & 38 & 1 & 10\hspace*{2pt}\\
& $>$25 & 88 & 89 & 0 & 6 & 6 & 12 & 23 & 29 & 6 & 22\hspace*{2pt}\\
\hline
\end{tabular*}\vspace*{-3pt}
\end{table}
identity. The behavior shown in
Tables~\ref{tabquartiels4}, \ref{tabquartiels44}, \ref{tabquartiels5} and~\ref{tabquartiels55}
is as in the previous two cases. The difference in the sparse model is
perhaps the most striking one.

\section{Proofs and ramification}\label{secyulewalkerdetails}

In this section, we will prove Theorems~\ref{thmyeconvrgextreme},
\ref{thmestdistrib},~\ref{thmexpandK}, and also explicitly mention
some auxiliary results which have interest in themselves.
%
%
\begin{table}[b]
\vspace*{-3pt}
\def\arraystretch{0.9}
\caption{Simulation of an $\operatorname{AR}(24)$ process with nonzero coefficients
$\theta_1 = 0.6$, $\theta_2 = -0.1$, $\theta_4 = 0.05$, $\theta_7 =
0.15$, $\theta_8 = -0.27$, $\theta_{10} = 0.1$, $\theta_{12} = -0.2$,
$\theta_{15} = -0.25$, $\theta_{18} = 0.05$, $\theta_{20} = 0.1$,
$\theta_{21} = -0.3$, $\theta_{24} = 0.17$, $\varepsilon\sim\Gaussian
(0,1)$, 1000 repetitions, $d_n \in\{38,42\}$}
\label{tabquartiels44}
\begin{tabular*}{\tablewidth}{@{\extracolsep{\fill}}lrrrrrrrrrrr@{}}
\hline
$\bolds{n}$& \multicolumn{1}{c}{$\widehat{\bolds{q}}$} &
\multicolumn{1}{c}{\textbf{AIC}\hspace*{-1pt}}
& \multicolumn{1}{c}{\textbf{AIC*}\hspace*{-3pt}}
& \multicolumn{1}{c}{\textbf{BIC}\hspace*{-1pt}}
& \multicolumn{1}{c}{\textbf{BIC*}\hspace*{-3pt}}
& \multicolumn{1}{c}{\textbf{HQC}\hspace*{-3pt}}
& \multicolumn{1}{c}{\textbf{HQC*}\hspace*{-5pt}}
& \multicolumn{1}{c}{\textbf{MIC}\hspace*{-2pt}}
& \multicolumn{1}{c}{\textbf{MIC*}\hspace*{-4pt}}
& \multicolumn{1}{c}{$\bolds{\widehat{q}_{y_n}^{(5)}}$\hspace*{-2pt}}
& \multicolumn{1}{c@{}}{$\bolds{\widehat{q}_{x_n}^{(5)}}$} \\
\hline
\hphantom{0}500 & $<$23 & 63 & 62 & 716 & 545 & 302 & 288 & 210 & 205 & 589 & 430\hspace*{2pt}\\
& 23 & 38 & 38 & 55 & 60 & 87 & 87 & 85 & 85 & 58 & 71\hspace*{2pt}\\
& \textbf{24} & \textbf{513} & \textbf{512} & \textbf{208} & \textbf{357} & \textbf{490} & \textbf{500} &
\textbf{525}
& \textbf{526} & \textbf{326} & \textbf{437}\hspace*{2pt}\\
& 25 & 192 & 192 & 18 & 28 & 93 & 93 & 129 & 129 & 19 & 27\hspace*{2pt}\\
& $>$25 & 194 & 196 & 3 & 10 & 28 & 32 & 51 & 55 & 8 & 35\hspace*{2pt}\\
[6pt]
1000 & $< $23 & 0 & 0 & 81 & 30 & 6 & 5 & 3 & 3 & 42 & 18\hspace*{2pt}\\
& 23 & 0 & 0 & 34 & 31 & 8 & 7 & 6 & 6 & 48 & 35\hspace*{2pt}\\
& \textbf{24} & \textbf{562} & \textbf{552} & \textbf{835} & \textbf{857} & \textbf{796} & \textbf{775} &
\textbf{761}
& \textbf{741} & \textbf{868} & \textbf{842}\hspace*{2pt}\\
& 25 & 197 & 195 & 48 & 45 & 140 & 137 & 160 & 156 & 7 & 24\hspace*{2pt}\\
& $>$25 & 241 & 253 & 2 & 37 & 50 & 76 & 70 & 94 & 35 & 81\hspace*{2pt}\\
\hline
\end{tabular*}
\end{table}
For ${d_n} \leq m$ let ${\bolds\Gamma}_m^{-1} = (\gamma_{i,j}^*
)_{1 \leq i,j \leq m}$ be the inverse of the covariance matrix ${\bolds
\Gamma}_m = (\gamma_{i,j} )_{1 \leq i,j \leq m}$ associated
to the $\operatorname{AR}(d_n)$ process~$\{X_k\}_{k \in\Z}$. Due to
Galbraith and Galbraith~\cite{galbraith1974}, it holds that
%
%
\begin{equation}\label{eqgammainversematrixexact}
\sigma^2 \gamma_{i,j}^* = \sum_{r = 0}^{\alpha} \theta_r \theta
_{r + j
- i} - \sum_{r = \beta}^{{d_n} + i - j} \theta_r \theta_{r + j -
i},\qquad
1 \leq i \leq j \leq m,
\end{equation}
where
\[
\alpha= \min\{i-1,{d_n}+i-j,m-j \},\qquad \beta= \max
\{
i-1,m-j \},
\]
and either of the sums is taken to be zero if its upper limit is less
than its lower
limit. The second sum is zero unless $m - d_n + 1 \leq i \leq j \leq d_n$
%
%
\begin{table}
\def\arraystretch{0.9}
\caption{Simulation of an $\operatorname{AR}(24)$ process with nonzero coefficients
$\theta_1 = 0.6$, $\theta_2 = -0.1$, $\theta_4 = 0.05$, $\theta
_{10} =
0.1$, $\theta_{12} = -0.2$, $\theta_{24} = 0.17$, $\varepsilon\sim
\Gaussian(0,1)$, 1000 repetitions, $d_n \in\{29,34\}$}
\label{tabquartiels5}
\begin{tabular*}{\tablewidth}{@{\extracolsep{\fill}}lrrrrrrrrrrr@{}}
\hline
$\bolds{n}$& \multicolumn{1}{c}{$\widehat{\bolds{q}}$} &
\multicolumn{1}{c}{\textbf{AIC}\hspace*{-1pt}}
& \multicolumn{1}{c}{\textbf{AIC*}\hspace*{-3pt}}
& \multicolumn{1}{c}{\textbf{BIC}\hspace*{-1pt}}
& \multicolumn{1}{c}{\textbf{BIC*}\hspace*{-3pt}}
& \multicolumn{1}{c}{\textbf{HQC}\hspace*{-3pt}}
& \multicolumn{1}{c}{\textbf{HQC*}\hspace*{-5pt}}
& \multicolumn{1}{c}{\textbf{MIC}\hspace*{-2pt}}
& \multicolumn{1}{c}{\textbf{MIC*}\hspace*{-4pt}}
& \multicolumn{1}{c}{$\bolds{\widehat{q}_{y_n}^{(5)}}$\hspace*{-2pt}}
& \multicolumn{1}{c@{}}{$\bolds{\widehat{q}_{x_n}^{(5)}}$} \\
\hline
125 & $<$23 & 1000 & 991 & 1000 & 991 & 1000 & 991 & 1000 & 991 & 991
& 969\hspace*{2pt}\\
& 23 & 0 & 2 & 0 & 2 & 0 & 2 & 0 & 2 & 2 & 6\hspace*{2pt}\\
&\textbf{24} & \textbf{0} & \textbf{6} & \textbf{0} & \textbf{6} & \textbf{0} & \textbf{6} & \textbf{0} & \textbf{6} &
\textbf{6} & \textbf{20}\hspace*{2pt}\\
& 25 & 0 & 1 & 0 & 1 & 0 & 1 & 0 & 1 & 1 & 1\hspace*{2pt}\\
& $>$25 & 0 & 0 & 0 & 0 & 0 & 0 & 0 & 0 & 0 & 4\hspace*{2pt}\\
[6pt]
250 & $<$23 & 857 & 768 & 1000 & 817 & 998 & 817 & 986 & 815 & 817 &
702\hspace*{2pt}\\
& 23 & 1 & 15 & 0 & 27 & 0 & 26 & 0 & 25 & 27 & 39\hspace*{2pt}\\
& \textbf{24} & \textbf{99} & \textbf{166} & \textbf{0} & \textbf{142} & \textbf{2} & \textbf{143} &
\textbf{13} & \textbf{145} & \textbf{142} & \textbf{225}\hspace*{2pt}\\
& 25 & 20 & 22 & 0 & 3 & 0 & 3 & 0 & 3 & 3 & 5\hspace*{2pt}\\
& $>$25 & 23 & 29 & 0 & 11 & 0 & 11 & 1 & 12 & 11 & 29\hspace*{2pt}\\
\hline
\end{tabular*}
\end{table}
while both sums are zero if $j - i > d_n$. Note that this implies
$\sigma^2(m) \gamma_{m,m}^* = 1$ for $m > d_n$, and in particular that
%
%
\begin{equation}\label{eqsupcovgoestozero}
{\sup_{|h|\geq n} \sup_{i} }|\gamma_{i,i + h}^* | =
\OOO
((\log n)^{-1} ),
\end{equation}
if Assumption~\ref{assmain} is valid.
Throughout this section and particularly in the proofs of the presented
%
%
\begin{table}[b]
\def\arraystretch{0.9}
\caption{Simulation of an $\operatorname{AR}(24)$ process with nonzero coefficients
$\theta_1 = 0.6$, $\theta_2 = -0.1$, $\theta_4 = 0.05$, $\theta
_{10} =
0.1$, $\theta_{12} = -0.2$,
$\theta_{24} = 0.17$, $\varepsilon\sim\Gaussian(0,1)$, 1000 repetitions,
$d_n \in\{38,42\}$}
\label{tabquartiels55}
\begin{tabular*}{\tablewidth}{@{\extracolsep{\fill}}lrrrrrrrrrrr@{}}
\hline
$\bolds{n}$& \multicolumn{1}{c}{$\widehat{\bolds{q}}$} &
\multicolumn{1}{c}{\textbf{AIC}\hspace*{-1pt}}
& \multicolumn{1}{c}{\textbf{AIC*}\hspace*{-3pt}}
& \multicolumn{1}{c}{\textbf{BIC}\hspace*{-1pt}}
& \multicolumn{1}{c}{\textbf{BIC*}\hspace*{-3pt}}
& \multicolumn{1}{c}{\textbf{HQC}\hspace*{-3pt}}
& \multicolumn{1}{c}{\textbf{HQC*}\hspace*{-5pt}}
& \multicolumn{1}{c}{\textbf{MIC}\hspace*{-2pt}}
& \multicolumn{1}{c}{\textbf{MIC*}\hspace*{-4pt}}
& \multicolumn{1}{c}{$\bolds{\widehat{q}_{y_n}^{(5)}}$\hspace*{-2pt}}
& \multicolumn{1}{c@{}}{$\bolds{\widehat{q}_{x_n}^{(5)}}$} \\
\hline
500 & $<$23 & 351 & 270 & 1000 & 383 & 952 & 380 & 854 & 379 & 383 &
256\hspace*{2pt}\\
& 23 & 2 & 8 & 0 & 51 & 0 & 48 & 0 & 41 & 51 & 61\hspace*{2pt}\\
& \textbf{24} & \textbf{451} & \textbf{522} & \textbf{0} & \textbf{547} & \textbf{45} & \textbf{550} &
\textbf{130} &
\textbf{545} & \textbf{547} & \textbf{637}\hspace*{2pt}\\
& 25 & 74 & 73 & 0 & 0 & 3 & 3 & 13 & 13 & 0 & 2\hspace*{2pt}\\
& $>$25 & 122 & 127 & 0 & 19 & 0 & 19 & 3 & 22 & 19 & 44\hspace*{2pt}\\
[6pt]
1000 & $< $23 & 10 & 6 & 986 & 15 & 440 & 15 & 280 & 15 & 15 & 3\hspace*{2pt}\\
& 23 & 0 & 0 & 0 & 14 & 0 & 13 & 0 & 11 & 14 & 12\hspace*{2pt}\\
& \textbf{24} & \textbf{718} & \textbf{715} & \textbf{14} & \textbf{941} & \textbf{522} & \textbf{908} &
\textbf{659}
& \textbf{887} & \textbf{941} & \textbf{905}\hspace*{2pt}\\
& 25 & 121 & 118 & 0 & 3 & 32 & 31 & 46 & 45 & 3 & 8\hspace*{2pt}\\
& $>$25 & 151 & 161 & 0 & 27 & 6 & 33 & 15 & 42 & 27 & 72\hspace*{2pt}\\
\hline
\end{tabular*}
\end{table}
results, we use the notation $\widehat{\sigma}^2 = \widehat{\sigma}^2(d_n)$.
Note that we can rewrite the equation defining the $\operatorname{AR}(d_n)$ process as
%
%
\begin{equation}
{\mathbf Y } = {\mathbf X} {\bolds\Phi_{d_n} } + {\mathbf Z},
\end{equation}
where ${\mathbf Y } = (X_1,\ldots,X_n)^T$, ${\mathbf Z} =
(\varepsilon
_1,\ldots,\varepsilon_n )^T$, and the
$n \times d_n$ design matrix~${\mathbf X}$ is given as
\[
{\mathbf X} =
\pmatrix{
X_0 & X_{-1} & \cdots& X_{1 - d_n}
\cr
X_1 & X_{0} & \cdots&
X_{2 - d_n} \cr
\cdots& \cdots& & &
\cr
X_{n-1} & X_{n-2} & \cdots& X_{n - d_n}}.
\]
We have
\[
\bolds\Gamma^{-1} \mathbf{X}^T \mathbf{Z} = { \bolds\Gamma}^{-1} \sum_{k = 1}^n
{\mathbf V}_k
= \sum_{k = 1}^n {\mathbf U}_k,
\]
where ${\mathbf V}_k = (V_{k}^{(1)},\ldots, V_{k}^{(d_n)})^T$,
${\mathbf
U}_k = (U_{k}^{(1)},\ldots, U_{k}^{(d_n)})^T$.
The following results are key ingredients.\vspace*{-3pt}
%
%
\begin{lem}\label{lemywinversemartrixtext}
Let $\{X_k\}_{k \in\Z}$ be an $\operatorname{AR}(d_n)$ process, such that
Assumption~\ref{assmain} is valid.
Then
\[
P\bigl(\|{\bolds\Gamma}_{d_n}^{-1} - {{\widehat{\bolds\Gamma}_{d_n}
}}^{-1}\|_{\infty} > (\log n)^{-\chi_1} \bigr) = \OO
\biggl(\frac
{(d_n( \log n)^{\chi_1})^{p}}{n^{p/2}} \biggr),\qquad \chi_1
\geq0.\vspace*{-3pt}\vadjust{\goodbreak}
\]
\end{lem}
%
%
\begin{lem}\label{lemapproxestimator}
Assume that the assumptions of Theorem~\ref{thmyeconvrgextreme} are valid.
Then we have
\[
P\Biggl( \Biggl\| n^{1/2}(\widehat{{\bolds\Theta}}_{d_n} - {\bolds\Theta
}_{d_n}) - n^{-1/2} \sum_{k = 1}^n {\mathbf U}_k \Biggr\|_{\infty} \geq
(\log n)^{-\chi_1} \Biggr) = \OOO
(1),
\]
where $1 < \chi_1$.
\end{lem}
%
%
\begin{lem}\label{lemwuspecialcase}
Assume that the assumptions of Theorem~\ref{thmyeconvrgextreme} are
valid. Then:
\begin{eqnarray*}
\mbox{\textup{(i)}}\quad
\lim_{n \to\infty} P\Biggl(\max_{1 \leq h \leq d_n}\sigma^{-1}
\Biggl|(n \gamma_{h,h}^*)^{-1/2}\sum_{k = 1}^n U_{k}^{(h)} \Biggr|
\leq u_n \Biggr) &=& \exp(-{\exp}(-x)),\\
\mbox{\textup{(ii)}\hspace*{135pt}}\quad
\sqrt{n}\|\widehat{\bolds\Phi}_{d_n} - {\bolds\Phi}_{d_n}
\|_{\infty} &=& \OO_P\bigl(\sqrt{\log d_n}\bigr),
\end{eqnarray*}
where $u_n = a_n z + b_n$, $a_n, b_n,z$ are as in
Theorem~\ref{thmyeconvrgextreme}.\vspace*{-2pt}
\end{lem}

The proofs of Lemmas~\ref{lemywinversemartrixtext},
\ref{lemapproxestimator} and~\ref{lemwuspecialcase}
are given in Section~\ref{secproofsofyw}. Based on the above results, one readily
derives the following weak version of Theorem~\ref{thmyeconvrgextreme}.\vspace*{-2pt}

%
\begin{cor}\label{corthmextremebasic}
Assume that the assumptions of Theorem~\ref{thmyeconvrgextreme} are valid.
Then for $z \in\R$
\[
P\Bigl(a_n^{-1}\Bigl({\sqrt{n}\max_{1 \leq h \leq d_n}} |
(\gamma_{h,h}^* \sigma^2)^{-1/2} (\widehat{\theta}_{i} - \theta
_{i})| - b_{n}\Bigr) \leq z \Bigr) \rightarrow\exp(- e^{-z}),
\]
where $a_{n}$ and $b_{n}$ are as in Theorem~\ref{thmyeconvrgextreme}.\vspace*{-2pt}
\end{cor}
%
%
\begin{cor}\label{corsigmaest}
Under the same conditions as in Theorem~\ref{thmyeconvrgextreme},
we have
\[
|\widehat{\sigma}^2 - \sigma^2| = \OO_P(n^{-1/2}
\log n
).\vspace*{-2pt}
\]
\end{cor}

Throughout the proofs, the following inequality will be frequently
used. For random variables $X_1,\ldots,X_q$, and $\varepsilon> 0$, the
inequality between the geometric and arithmetic mean implies
%
%
\begin{equation}\label{eqinequ}
P\Biggl( \prod_{i = 1}^q |X_i| \geq\varepsilon\Biggr) \leq\sum_{i = 1}^q
P( |X_i| \geq\varepsilon^{1/q}).\vspace*{-2pt}
\end{equation}
\begin{pf*}{Proof of Corollary~\ref{corthmextremebasic}}
It holds that
\begin{eqnarray*}
&& \max_{1 \leq h \leq d_n}\sigma^{-1} \Biggl|(n \gamma
_{h,h}^*)^{-1/2}\Biggl(n (\widehat{\theta}_h - \theta_h) - \sum_{k =
1}^n U_{k}^{(h)} \Biggr) \Biggr|\\[-2pt]
&&\qquad \leq\sqrt{\Bigl(n \sigma^2 \inf_{h}
\gamma
_{h,h}^*\Bigr)^{-1}} \max_{1 \leq h \leq d_n}\Biggl|\Biggl(n (\widehat
{\theta
}_h - \theta_h) - \sum_{k = 1}^n U_{k}^{(h)} \Biggr) \Biggr|.
\end{eqnarray*}
Since $\inf_{h} \gamma_{h,h}^* > 0$ and choosing $\chi_1 > 1$, the
claim follows from Lemmas~\ref{lemapproxestimator}
and~\ref{lemwuspecialcase}.\vspace*{-2pt}
\end{pf*}
\begin{pf*}{Proof of Corollary~\ref{corsigmaest}}
Trivially, it holds that
\begin{eqnarray*}
\widehat{\sigma}^2 - \sigma^2 &=& \widehat{\phi}_0 - \phi_0 +
\widehat
{\bolds\Theta}_{d_n}^T \widehat{\bolds\Phi}_{d_n} - {\bolds\Theta}_{d_n}^T
{\bolds\Phi}_{d_n}\\[-2pt]
&=& \widehat{\phi}_0 - \phi_0 + (\widehat
{\bolds
\Theta}_{d_n}^T - {\bolds\Theta}_{d_n}^T)(\widehat{\bolds
\Phi
}_{d_n} - {\bolds\Phi}_{d_n})\\[-2pt]
&&{} + (\widehat{\bolds\Theta
}_{d_n}^T -
{\bolds\Theta}_{d_n}^T){\bolds\Phi}_{d_n} + {\bolds\Theta
}_{d_n}^T
(\widehat{\bolds\Phi}_{d_n} - {\bolds\Phi}_{d_n}).
\end{eqnarray*}

By Corollary~\ref{corthmextremebasic} and Lemma~\ref{lemwuspecialcase} we have
\begin{eqnarray*}
\|(\widehat{\bolds\Theta}_{d_n}^T - {\bolds\Theta
}_{d_n}^T
)(\widehat{\bolds\Phi}_{d_n} - {\bolds\Phi}_{d_n}) \|
_{\infty
} &\leq& d_n\|\widehat{\bolds\Theta}_{d_n}^T - {\bolds\Theta
}_{d_n}^T\|
_{\infty} \|\widehat{\bolds\Phi}_{d_n} - {\bolds\Phi}_{d_n}
\|
_{\infty}\\[-2pt]
&=& \OO_P(n^{-1/2} \log n).\vadjust{\goodbreak}
\end{eqnarray*}

Similarly, we obtain from Lemmas~\ref{lemwuspecialcase},~\ref{lemboundmatrixnorms} and Assumption~\ref{assmain}
\begin{eqnarray*}
\|(\widehat{\bolds\Theta}_{d_n}^T - {\bolds\Theta
}_{d_n}^T
){\bolds\Phi}_{d_n}\|_{\infty} &=& \OO_P(n^{-1/2} \log n
),\\
\|{\bolds\Theta}_{d_n}^T(\widehat{\bolds\Phi}_{d_n} -
{\bolds
\Phi}_{d_n})\|_{\infty} &=& \OO_P(n^{-1/2} \log n
).
\end{eqnarray*}
Moreover, from the above one readily deduces $|\widehat{\phi}_0 -
\phi_0 | = \OO_P(n^{-1/2} \log n )$. Piecing everything
together, the claim follows.
\end{pf*}
\begin{pf*}{Proof of Theorem~\ref{thmyeconvrgextreme}}
Due to Corollary~\ref{corthmextremebasic}, it suffices to show that
the error difference
%
%
\begin{eqnarray}\quad
\max_{1 \leq i \leq d_n} \Delta_i &=& \max_{1 \leq i \leq d_n}\sqrt
{n}|(\widehat{\gamma}_{i,i}^* \widehat{\sigma
}^2)^{-1/2}(\widehat
{\theta}_{i} - \theta_{i}) - (\gamma_{i,i}^* \sigma
^2)^{-1/2}(\widehat
{\theta}_{i} - \theta_{i})| \nonumber\\[-8pt]\\[-8pt]
&=& \OO_P( (\log n)^{-\chi
_1})\nonumber
\end{eqnarray}
for some $\chi_1 > 1$.
Note that per assumption we have that $(\log n)^{\chi_2 p} n^{-p/2}
d_n^p = \OOO(1)$ for some
$\chi_2 > 1$.
Moreover,
\begin{eqnarray*}
\max_{0 \leq i \leq d_n} \Delta_i &\leq&\max_{0 \leq i \leq d_n}
\bigl|
\bigl((\widehat{\gamma}_{i,i}^* \widehat{\sigma}^2)^{1/2} -
(\gamma
_{i,i}^* \sigma^2)^{1/2} \bigr) (\widehat{\gamma}_{i,i}^* \widehat
{\sigma}^2 )^{-1/2}\bigr| \\
&&{}\times{\sqrt{n}\max_{0 \leq i \leq d_n}}
|(\widehat{\theta}_{i} - \theta_{i})(\gamma_{i,i}^* \sigma
^2)^{-1/2}|.
\end{eqnarray*}
Corollary~\ref{corthmextremebasic} gives us ${\sqrt{n}\max_{0 \leq i
\leq d_n}} |(\widehat{\theta}_{i} - \theta_{i})(\gamma_{i,i}^*
\sigma^2)^{-1/2}| = \OO_P(\log n )$, hence we need to
study $|(\widehat{\gamma}_{i,i}^* \widehat{\sigma}^2)^{1/2} -
(\gamma_{i,i}^* \sigma^2)^{1/2} | (\widehat{\gamma}_{i,i}^*
\widehat{\sigma}^2 )^{-1/2}$. Since
\begin{eqnarray*}
\biggl|\frac{(\widehat{\gamma}_{i,i}^* \widehat{\sigma}^2)^{1/2} -
(\gamma_{i,i}^* \sigma^2)^{1/2}}{(\widehat{\gamma}_{i,i}^* \widehat
{\sigma}^2 )^{1/2}}\biggr|
&=& \biggl|\frac{\widehat{\gamma}_{i,i}^*
\widehat{\sigma}^2 - \gamma_{i,i}^* \sigma^2}{(\widehat{\gamma}_{i,i}^*
\widehat{\sigma}^2)^{1/2} + (\gamma_{i,i}^* \sigma^2)^{1/2}} \frac
{1}{(\widehat{\gamma}_{i,i}^* \widehat{\sigma}^2 )^{1/2}}\biggr|\\
&\leq&
\biggl|\frac{(\widehat{\gamma}_{i,i}^* \widehat{\sigma}^2 -
\gamma
_{i,i}^* \sigma^2)}{\widehat{\gamma}_{i,i}^* \widehat{\sigma}^2}
\biggr|,
\end{eqnarray*}
it suffices to treat $(\widehat{\gamma}_{i,i}^* \widehat
{\sigma}^2
- \gamma_{i,i}^* \sigma^2)(\widehat{\gamma}_{i,i}^* \widehat
{\sigma
}^2 )^{-1}$. For $\varepsilon= \log n^{-\chi_2}$ we have
\begin{eqnarray*}
&&\{|\sigma^2 \gamma_{i,i}^* - \widehat{\sigma}^2\widehat
{\gamma
}_{i,i}^* | \geq\varepsilon\widehat{\gamma}_{i,i}^* \widehat
{\sigma
}^2 \} \\
&&\qquad\subseteq \{|\sigma^2 \gamma_{i,i}^* -
\widehat
{\sigma}^2\widehat{\gamma}_{i,i}^* |(1 + \varepsilon) \geq
\varepsilon
\gamma_{i,i}^* \sigma^2 \}\\
&&\qquad\subseteq \{|\sigma^2
\gamma
_{i,i}^* - \widehat{\sigma}^2\widehat{\gamma}_{i,i}^* | \geq
\varepsilon\gamma_{i,i}^* \sigma^2/2 \}\\
&&\qquad\subseteq \{
|\sigma^2 \widehat{\gamma}_{i,i}^* - \widehat{\sigma}^2\widehat
{\gamma
}_{i,i}^* | + |\sigma^2 \gamma_{i,i}^* - \sigma^2
\widehat
{\gamma}_{i,i}^* | \geq\varepsilon\gamma_{i,i}^* \sigma
^2/2\}.
\end{eqnarray*}
Since $\sigma^2, \gamma_{i,i}^* \geq C > 0$, we have from Lemma
\ref{lemywinversemartrixtext} that for $1 < \chi_1 < \chi_2$
\[
P\Bigl({\max_{0 \leq i \leq d_n}}|\sigma^2 \gamma_{i,i}^* -
\sigma
^2\widehat{\gamma}_{i,i}^* | \geq\log n^{-\chi_1} \min_{1
\leq i
\leq d_n}\gamma_{i,i}^*\Bigr) = \OO(\log n^{\chi_2 p} n^{-p/2}
d_n^p).
\]
In order to treat $|\sigma^2 \widehat{\gamma}_{i,i}^* -
\widehat
{\sigma}^2\widehat{\gamma}_{i,i}^* |$, note that
\[
{\max_{0 \leq i \leq d_n}} |\sigma^2 \widehat{\gamma}_{i,i}^* -
\widehat{\gamma}_{i,i}^* \widehat{\sigma}^2| \leq\max_{0
\leq i
\leq d_n} (|\sigma^2 - \widehat{\sigma}^2|
|\gamma
_{i,i}^* - \widehat{\gamma}_{i,i}^* | + \gamma_{i,i}^*
|\sigma
^2 - \widehat{\sigma}^2|),
\]
which by virtue of Corollary~\ref{corsigmaest} and Lemma~\ref{lemywinversemartrixtext} is of the magnitude
$\OO_P(n^{-1/2}\times
\log n )$. We thus obtain that
\[
P\Bigl(\max_{0 \leq i \leq d_n} \Delta_i \geq\log n^{-\chi_1}
\Bigr) =
\OOO(1)
\]
for some $\chi_1 > 1$, which completes the proof.
\end{pf*}
\begin{pf*}{Proof of Corollary~\ref{corsimple}}
First note that both conditions (i) and (ii) imply that $|\alpha_i| =
\OO(\rho^{-i})$, $0 < \rho< 1$ (cf.
\cite{timeseriesbrockwell}). Hence Remark~\ref{remdiscusscondis} yields
that we may choose $d_n = \OO(n^{\delta})$, $0< \delta<
1/2$. Now assume that (i) holds. Then relation~(\ref{eqgammainversematrixexact}) implies
\[
\sigma^2 \inf_h \gamma_{h,h}^* \geq1 - \sum_{i = 1}^{d_n} |\theta_i|^2
\geq1 - \sum_{i = 1}^{d_n} |\theta_i| > 0,
\]
whence the claim. If (ii) holds, then for large enough $n$ we obtain similarly
\[
\sigma^2 \inf_h \gamma_{h,h}^* \geq\sum_{i = 0}^{\alpha} \theta
_i^2 -
\sum_{i = \beta}^{d_n} \theta_i^2 \geq\sum_{i = 0}^{\alpha}
\theta_i^2
\geq1,
\]
where $\alpha$, $\beta$ are as in~(\ref{eqgammainversematrixexact}).
\end{pf*}

We are now ready to prove Theorem~\ref{thmestdistrib}.
\begin{pf*}{Proof of Theorem~\ref{thmestdistrib}}
Let $q_0 = q$ be the true order of the $\operatorname{AR}(q)$-process $\{
X_k
\}_{k \in\Z}$, put
\[
\overline{\theta}_{i,n} = a_n^{-1}\bigl(\sqrt{n} |(\widehat
{\gamma
}_{i,i}^* \widehat{\sigma}^2)^{-1/2}(\widehat{\theta}_{i} - \theta
_{i})| - b_{n}\bigr)
\]
and assume first that $k \in\N$, $k > 0$. Note that $\theta_i = 0$ for
$i > q$. Then we have that
\begin{eqnarray*}
P(\widehat{q}_{z_n} = k + q ) &=& P\Bigl(\{\overline
{\theta
}_{q + k,n} > z_n \} \cap\Bigl\{\max_{k + q + 1 \leq i \leq d_n}\overline
{\theta}_{i,n} \leq z_n\Bigr\} \Bigr)\\
&=& P\Bigl(\max_{k + q \leq i \leq
d_n}\overline{\theta}_{i,n} \leq z_n \Bigr) - P\Bigl(\max_{k + q + 1
\leq i \leq d_n}\overline{\theta}_{i,n} \leq z_n \Bigr).
\end{eqnarray*}

Due to Theorem~\ref{thmgaussianapprox}, we can approximate the sequence
$\{\overline{\theta}_{i,n} \}_{1 \leq i \leq d_n}$ by a~suitably
transformed corresponding\vspace*{1pt} sequence of mean-zero Gaussian
random variables $\xiv_{d_n} = (\xi_{n,1},\ldots, \xi_{n,d_n} )^T$ with
covariance matrix ${\bolds\Gamma}_{\xiv_{d_n}}^*$. Let\vspace*{-1pt}
$\etav _{d_n} = (\eta_{n,1},\ldots, \eta_{n,d_n} )^T$ be another
sequence of i.i.d. mean-zero Gaussian random variables with unit
variance. Following Deo~\cite{deoabs}, we obtain from
${\max}|{\bolds\Gamma }_{\xiv_{d_n}}^* - {\bolds\Gamma}_{d_n}^* | =
\OOO(d_n^{-1} )$ that for fixed $l \in\N$
\begin{eqnarray*}
&&\Bigl|P\Bigl(\max_{q + l\leq i \leq d_n}a_n^{-1}(|\xi_{n,i}
| - b_{n}) \leq z_n \Bigr) - P\Bigl(\max_{q + l\leq i
\leq
d_n}a_n^{-1}(|\eta_{n,i} | - b_{n}) \leq z_n
\Bigr)\Bigr| \\
&&\qquad \leq
C \sum_{1 \leq i < j \leq d_n} |\rho_{i,j} |
\bigl(d_n^{-
{2 z_n^2}/({1 + |\rho_{i,j}|})}\bigr).
\end{eqnarray*}
Imitating the technique in Berman~\cite{Berman1964}, we obtain that
the above quantity is of the magnitude $\OOO(d_n^{(-z_n^2 + 1)/2}
)$. This yields
\begin{eqnarray*}
&&
P(\widehat{q}_{z_n} = k + q ) \\
&&\qquad= P\Bigl(\max_{q +
k +
1\leq i \leq d_n}a_n^{-1}(|\eta_{n,i} | - b_{n})
\leq z_n \Bigr)\\
&&\qquad\quad{} - P\Bigl(\max_{q + k \leq i \leq d_n}a_n^{-1}
(
|\eta_{n,i} | - b_{n}) \leq z_n \Bigr) +
\OOO\bigl(n^{-\nu} +
d_n^{(-z_n^2 + 1)/2} \bigr)\\
&&\qquad= P\bigl(a_n^{-1}(|\eta_{n,1}
| - b_{n})\leq z_n \bigr)^{d_n - k - q}\bigl(1 - P
\bigl(a_n^{-1}(|\eta_{n,1} | - b_{n})\leq z_n
\bigr)\bigr) \\
&&\qquad\quad{} + \OOO\bigl(n^{-\nu} +
d_n^{(-z_n^2 + 1)/2} \bigr).
\end{eqnarray*}
From the definition of $a_n, b_n$, and since $z_n \to\infty$, we
obtain that (Deo~\cite{deoabs})
%
%
\begin{eqnarray}
&\displaystyle \lim_n P\bigl(a_n^{-1}(|\eta_{n,1} | - b_{n}
)\leq
z_n \bigr)^{d_n - k - q} \to1,& \\
&\displaystyle P\bigl(a_n^{-1}(|\eta_{n,1} | - b_{n}) > z_n
\bigr) = \frac{e^{-z_n}}{d_n} + \OOO
\biggl(\frac{e^{-z_n}}{d_n}\biggr).&
\end{eqnarray}
This yields
%
%
\begin{equation}
P(\widehat{q}_{z_n} = k + q ) = \frac{e^{-z_n}}{d_n} +
\OOO
\biggl(\frac{e^{-z_n}}{d_n} + d_n^{(-z_n^2 + 1)/2}\biggr)
\end{equation}
and in particular
%
%
\begin{equation}\label{eqconsistencyproof}
P(\widehat{q}_{z_n} > q ) = \sum_{k = 1}^{d_n} P
(\widehat
{q}_{z_n} = k + q ) = e^{-z_n} +
\OOO(e^{-z_n} + d_n^{-z_n^2 +
2}),
\end{equation}
and per assumption the right-hand side goes to zero as $n$ increases.
We now consider the case $P(\widehat{q}_{z_n} < q )$. To this
end, let $k \in\N$, $k > 0$. Then we have
\begin{eqnarray*}
P(\widehat{q}_{z_n} = q - k ) &\leq& P(\overline
{\theta
}_{q - k,n} \leq z_n ) \\
&=& P\bigl(a_n^{-1}\bigl(\bigl|\xi_{n,q
- k}
+ \sqrt{n}\theta_{q - k}\bigr| - b_{n}\bigr)\leq z_n \bigr) \\
&&{}+ \OO
(n^{-\nu}).
\end{eqnarray*}
Since $|\theta_{q - k}| > 0$, one readily verifies by known properties
of the Gaussian c.d.f. that $P(a_n^{-1}(|\xi_{n,q - k} +
\sqrt{n}\theta_{q - k}| - b_{n})\leq z_n ) = \OO
(n^{-\nu})$, and hence
%
%
\begin{equation}
P(\widehat{q}_{z_n} = q - k ) = \OO(n^{-\nu})
\end{equation}
and in particular
%
%
\begin{equation}
P(\widehat{q}_{z_n} < q ) = \OO(d_n n^{-\nu}
) \to0
\end{equation}
as $n$ increases. This together with~(\ref{eqconsistencyproof})
establishes consistency.
\end{pf*}
\begin{pf*}{Proof of Theorem~\ref{thmexpandK}}
Let $q_0 = q$ be the true order of the $\operatorname{AR}(q)$-process
$\{X_k \} _{k \in\Z}$. The proof then consists of two parts. It is
first shown that $P (\widehat{q}_n^* < q) \to0$, whereas in the second
part the claim $P(\widehat{q}_n^* > q) \to0$ is established.

First note that Lemma~\ref{lemboundmatrixnorms} and the Cauchy
interlacing theorem yield that $\|{\bolds\Gamma}_{k} \|
_{\infty
}$, $\|{\bolds\Gamma}_{k}^{-1}\|_{\infty} \leq C < \infty$,
uniformly for $1 \leq k \leq d_n$.
Hence, using that $\widehat{\bolds\Gamma}_{k}^{-1} \widehat{\bolds\Phi
}_k =
\widehat{\bolds\Theta}_k$, Lemma~\ref{lemwuspecialcase} and a~slight
adaption of Lemma~\ref{lemywinversemartrixtext} imply that
\[
|\widehat{\sigma}^2(k) - \sigma^2(k)| =
\OOO_P(1)\qquad
\mbox{uniformly for $1 \leq k \leq q$.}
\]
Since $\inf_h|1 - \theta_h^2| > 0$, we conclude that $\inf_k
\sigma^2(k) > 0$ and hence
\begin{equation}\label{eqbefhannan}
|{\log}(\widehat{\sigma}^2(k)) - \log(\sigma ^2(k) ) | =
\OOO_p(1).
\end{equation}
By Hannan~\cite{hannantimeseries}, Chapter VI, it holds that for $k
\in\N$
\begin{equation}\label{eqhannan}
\log(\widehat{\sigma}^2(k)) = \log\widehat{\phi}_{n,0} + \sum _{j =
1}^k \log\bigl(1 - \widehat{\theta}_j^2(k)\bigr).
\end{equation}
Then, arguing as in Hannan and Quinn~\cite{hannanarp}, we have due
to $C_n = \OOO(n)$ that for large enough $n$
\[
f_n(k) = \log(\widehat{\sigma}^2(k)) + n^{-1} C_n k
\]
is a~decreasing function in $k$ for $0 \leq k < q$, and strictly
decreasing for $q - 1 \leq k \leq q$ (since $\theta_q^2 > 0$) with
probability approaching one. This implies that eventually $\widehat
{q}_n^* \geq q$, hence it suffices to establish that the probability of
overestimating the order goes to zero as $n$ increases, that is,
\begin{equation}
\lim_n P\Bigl( \argmin_{q \leq k \leq d_n}\bigl(\log
(\widehat
{\sigma}^2(k)) + n^{-1} C_n k \bigr) \geq q + 1 \Bigr) = 0.
\end{equation}
Using the same arguments as in~\cite{hannan1982}, it follows that it
suffices to establish
\begin{equation}\label{eqthmorder1}\qquad
\lim_n P\Biggl(\max_{1 \leq k \leq d_n - q} \Biggl(\sum_{j = 1 +
q}^{k +
q} -\log\bigl(1 - \widehat{\theta}_j^2(k)\bigr) - n^{-1} C_n k \Biggr) \geq0
\Biggr) = 0.
\end{equation}
By Theorem~\ref{thmyeconvrgextreme}, we have that
\begin{equation}
\|\widehat{\bolds\Theta}_k^2\|_{\infty} = \OO_P(n^{-1}
\log d_n)\qquad \mbox{for $q_0 < k \leq d_n$.}
\end{equation}
This implies that for some increasing $\chi_n \to\infty$, we obtain that
\begin{equation}
-\sum_{j = 1 + q}^{k + q} \log\bigl(1 - \widehat{\theta}_j^2(k)\bigr) \leq k
\chi
_n n^{-1} \log d_n,
\end{equation}
with probability approaching one. Since $\log d_n = \OOO(C_n)$
per assumption, (\ref{eqthmorder1}) follows, which completes the proof.
\end{pf*}

\section{\texorpdfstring{Proofs of the auxiliary results of Section \protect\ref{secyulewalkerdetails}}{Proofs of the auxiliary results of Section 4}}\label{secproofsofyw}

The following result is required for the proofs.
%
%
\begin{lem}\label{lemproparp}
Let $\{X_k\}_{k \in\Z}$ be an $\operatorname{AR}(q)$ process such that
Assumption~\ref{assmain} is satisfied. Then:
\begin{longlist}
\item ${\sum_{h = 0}^{\infty}} |{\Cov}(X_k,X_{k + h}
)| < \infty$,\vspace*{2pt}
\item
$\sqrt{n} \|\widehat{\phi}_{n,h} - \phi_h\|
_p =
\OO(1)$, $p \geq1$.
\end{longlist}
\end{lem}
\begin{pf}
Both properties (i), (ii) follow from Assumption~\ref{assmain} via
straightforward computations (cf.
\cite{timeseriesbrockwell,hannantimeseries}).
\end{pf}

Recall the notation ${\bolds\Gamma}_m = (\gamma_{i,j} )_{1
\leq
i,j \leq m}$ and ${\bolds\Gamma}_m^{-1} = (\gamma_{i,j}^* )_{1
\leq i,j \leq m}$
for the covariance matrix and its inverse.
%
%
\begin{lem}\label{lemboundmatrixnorms}
Assume that Assumption~\ref{assmain} holds. Then for $d_n \leq m$ we have
$\|{\bolds\Gamma}_{m} \|_{\infty}$, $\|{\bolds\Gamma
}_{m}^{-1}\|_{\infty} \leq C < \infty$, uniformly in $m$.
\end{lem}
\begin{pf}
Using relation~(\ref{eqgammainversematrixexact}) and the
corresponding notation, one obtains
\begin{eqnarray*}
\|{\bolds\Gamma}_{m}^{-1}\|_{\infty} &=& \sigma^{-2} \max_{1
\leq j \leq m} \sum_{i = 1}^{m} |\sigma^2 \gamma_{i,j}^*| \\
&\leq&2
\sigma
^{-2} \max_{1 \leq j \leq m} \Biggl|\sum_{r = 0}^{\alpha} \theta_r
\theta_{r + j - i} - \sum_{r = \beta}^{{d_n} + i - j} \theta_r
\theta
_{r + j - i}\Biggr|\\&\leq&4 \sigma^{-2} \sum_{|h| \leq m} \sum_{r =
0}^{m} \bigl| \theta_r \theta_{|r + h|}\bigr| \leq8 \sigma^{-2}
\Biggl(\sum_{r = 0}^{\infty} |\theta_r| \Biggr)^2,
\end{eqnarray*}
where $\theta_h = 0$ for $h < 0$. Due to Assumption~\ref{assmain}, the
above expression is finite, hence the first claim follows. In order
to\vadjust{\goodbreak}
establish the result for ${\bolds\Gamma}_{m}$, note that
\[
\|{\bolds\Gamma}_{m}\|_{\infty}=\max_{1 \leq j \leq m}
\sum
_{i = 1}^{m} | \gamma_{i,j}| \leq2 \sum_{h = 0}^{\infty} |\phi_h| <
\infty
\]
by Lemma~\ref{lemproparp}(i), which yields the claim.
\end{pf}

We can now prove Lemma~\ref{lemwuspecialcase}, which we reformulate
below for the sake of readability.
%
%
\begin{lem}\label{lemwuapprox}
Suppose that $\inf_{h}|\gamma_{h,h}^* | > 0$ and Assumption
\ref{assmain} holds. Then:
\begin{eqnarray*}
\mbox{\textup{(i)}}\quad \lim_{n \to\infty} P\Biggl(\max_{0 \leq h \leq d_n}\sigma^{-1}
\Biggl|(n \gamma_{h,h}^*)^{-1/2}\sum_{k = 1}^n U_{k}^{(h)} \Biggr|
\leq u_n \Biggr)& = &\exp(-{\exp}(-x)),\\
\mbox{\textup{(ii)}}\quad\hspace*{135pt}
\sqrt{n}\|\widehat{\bolds\Phi}_{d_n} - {\bolds\Phi}_{d_n}
\|_{\infty} &=& \OO_P\bigl(\sqrt{\log d_n}\bigr).
\end{eqnarray*}
\end{lem}
\begin{pf}
We will first show (i). Using the notation established in Section~\ref{secyulewalkerdetails}, we have
\[
U_{{k}}^{(h)} = \varepsilon_{k} \Biggl(\sum_{j = 1}^{d_n} \gamma
_{h,j}^{(*)} \sum_{i = 0}^{\infty} \alpha_i \varepsilon_{k - j - i}
\Biggr) := \varepsilon_{k} \sum_{r = 1}^{\infty} \alpha_{r,h}^* \varepsilon_{k
- r},
\]
where $\alpha_{r,h}^* = \sum_{\{i \geq0, j \geq0, i +
j = r\}}\gamma_{h,j}^* \alpha_i$.
Let $0 < \delta< \delta^*$, and put $m_n = \lfloor n^{\delta^*}
\rfloor
$. Then it follows from Lemma~\ref{lemboundmatrixnorms} that
\[
\sup_{h} \sum_{r = m_n}^{\infty} |\alpha_{r,h}^*| \leq
C \sum
_{i = m_n - d_n}^{\infty} |\alpha_i| = \OO\bigl((m_n -
d_n)^{-\vartheta} \bigr) = \OO(m_n^{-\vartheta} ).
\]
Due to Assumption~\ref{assmain}, one may thus repeat the (quite
lengthy) proof of Theorem~1 (see also Remark 2) in
\cite{wuxiaobeats2011} to obtain the result. In fact, the present case is
easier to handle, since $\{U_{{k}}^{(h)}\}_{k \in\N}$ is a~martingale sequence.

Assertion (ii) follows directly from Theorem 1
in~\cite{wuxiaobeats2011}.
\end{pf}

We can now proof Lemma~\ref{lemywinversemartrixtext}, which we
restate for the sake of readability.

%
\begin{lem}\label{lemmatrixinvapproxerror}
If Assumption~\ref{assmain} holds, we have for $\chi_1 > 0$
\[
P\bigl(\|\widehat{\bolds\Gamma}_{d_n}^{-1} - {\bolds\Gamma}_{d_n}^{-1}
\|_{\infty} \geq(\log n)^{-\chi_1} \bigr) = \OO\biggl(\frac{(d_n
(\log n)^{\chi_1})^{p}}{n^{p/2}} \biggr).
\]
\end{lem}
\begin{pf}
We introduce the following abbreviations. Put
\[
E= \|{\bolds\Gamma}_{d_n}^{-1} \|_{\infty},\qquad F = \|
\widehat{\bolds\Gamma}_{d_n}^{-1} - {\bolds\Gamma}_{d_n}^{-1}\|
_{\infty
},\qquad G = \|\widehat{\bolds\Gamma}_{d_n} - {\bolds\Gamma}_{d_n}
\|
_{\infty}.
\]
Due to the stationarity of $\{X_k\}_{k \in\Z}$ it follows that
%
%
\begin{equation}\label{eqspecpropmarkov}
G= \|\widehat{\bolds\Gamma}_{d_n} - {\bolds\Gamma}_{d_n}\|
_{\infty
} \leq2 \sum_{h \leq d_n} \bigl|\widehat{\phi}_{n,|h|} - \phi
_{|h|}\bigr|,\vadjust{\goodbreak}
\end{equation}
and thus an application of the H\"{o}lder and Minikowski inequalities yields
%
%
\begin{equation}
\E(|G|) \leq2 \sum_{h \leq d_n} \bigl\|\widehat{\phi
}_{n,|h|} - \phi_{|h|}\bigr\|_p.
\end{equation}

Due to Lemma~\ref{lemproparp}(ii) we have $\sqrt{n}\|\widehat
{\phi}_{n,|i - j|} - \phi_{|i-j|} \|_p \leq C_p$ for some finite
constant $C_p$, thus the Markov inequality in connection with
Minikowski's inequality implies
%
%
\begin{equation}\label{eqspectwo}
P\bigl(\|\widehat{\bolds\Gamma}_{d_n} - {\bolds\Gamma}_{d_n}
\|
_{\infty} \geq(\log n)^{-\chi_1} \bigr) = \OO\biggl(\frac{(d_n
(\log
n)^{\chi_1})^{p}}{n^{p/2}} \biggr).
\end{equation}

Due to the sub-multiplicativity of the matrix norm \mbox{$\|\cdot\|_{\infty}$},
proceeding as in Lem\-ma 3 in~\cite{berk1978} one obtains
\[
F \leq(E+F) G E,
\]
and in particular if $E G < 1$
\[
F\leq E^2 G/(1 - E G).
\]
Since we have $E < \infty$ due to Lemma~\ref{lemboundmatrixnorms},
we deduce that for sufficiently large $n$
\[
P(F \geq\varepsilon) \leq P\bigl(G \geq(\log n)^{-1} \bigr)
+ P(G \geq E^2/2 \varepsilon).
\]
Choosing $\varepsilon= (\log n)^{-\chi_1}$, the claim follows.
\end{pf}

We are now in the position to show Lemma~\ref{lemywinversemartrixtext}. Recall that we have
%
%
\begin{equation}
{\mathbf Y } = {\mathbf X} {\bolds\Phi}_d + {\mathbf Z},
\end{equation}
where ${\mathbf Y } = (X_1,\ldots,X_n)^T$, ${\mathbf Z} =
(\varepsilon
_1,\ldots,\varepsilon_n )^T$, and ${\mathbf X}$ is the $n \times d_n$
design matrix.

We introduce the estimator $\widetilde{\bolds\Theta} =
(\widetilde
{\theta}_1,\ldots,\widetilde{\theta}_d)^T$ via
%
%
\begin{equation}
\widetilde{\bolds\Theta} = ({ \mathbf X}^T { \mathbf X})^{-1} { \mathbf
X}^T {\mathbf Y}.
\end{equation}

%
\begin{rem}\label{remmatrixinv}
It is evident from the proof that Lemma~\ref{lemmatrixinvapproxerror} remains valid if one replaces $\widehat
{\bolds\Gamma}_{d_n}$ with $n({\mathbf X}^T{\mathbf X})^{-1}$, which in
fact is the
better estimator.
\end{rem}
%
%
\begin{prop}\label{propapproxestimatorreally}
Let $\{X_k\}_{k \in\Z}$ be an $\operatorname{AR}(d_n)$ process, such that
the assumptions of Theorem~\ref{thmyeconvrgextreme} are satisfied. Then
\[
P\bigl(\bigl\| \sqrt{n}(\widehat{\bolds\Theta} -
\widetilde{\bolds
\Theta}) \bigr\|_{\infty} \geq(\log n)^{-\chi_1} \bigr) = \OO
((\log
n)^{\chi_1 p/2} n^{-p/4} d_n^{p/4 + 1}) + \OOO(1).
\]
\end{prop}
\begin{pf}
Following the proof of~\cite{timeseriesbrockwell}, Theorem 8.10.1, we
have the following decomposition:
\[
\sqrt{n}(\widehat{\bolds\Theta} - \widetilde{\bolds\Theta}
) =
\sqrt{n}{{\widehat{\bolds\Gamma}_{d_n} }}^{-1} ({\widehat
{\bolds\Phi
}_{d_n}} - n^{-1} {\mathbf X}^T{\mathbf Y}) + n^{1/2}\bigl({{\widehat
{\bolds\Gamma
}_{d_n} }}^{-1} - n({\mathbf X}^T{\mathbf X})^{-1} \bigr) n^{-1} {\mathbf
X}^T{\mathbf Y}.\vadjust{\goodbreak}
\]

For the $i$th component of $\sqrt{n} ({\widehat{\bolds\Phi
}_{d_n}} -
n^{-1} {\mathbf X}^T{\mathbf Y})$, which we denote with $\Upsilon_i$,
we have
\[
n^{-1/2} \sum_{k = 1-i}^0 X_k X_{k + i} + \sqrt{n} \overline{X}_n
\Biggl( (1 - n^{-1} i ) \overline{X}_n - n^{-1} \sum_{k = 1}^{n - i} (X_k +
X_{k + i}) \Biggr).
\]
Using the Minikowski and the Cauchy--Schwarz inequalities we get
\begin{eqnarray*}
&&\Biggl\|n^{-1/2} \sum_{k = 1-i}^0 X_k X_{k + i} + \sqrt{n} \overline
{X}_n \Biggl( (1 - n^{-1} i ) \overline{X}_n - n^{-1} \sum_{k = 1}^{n -
i} (X_k + X_{k + i}) \Biggr) \Biggr\|_{p/2} \\
&&\qquad\leq\sqrt{\frac{|1 -
i|}{n}}\Biggl\||1-i|^{-1/2} \sum_{k = 1-i}^0 (X_k X_{k + i} - \phi_{i})
\Biggr\|_{p/2} + n^{-1/2} \sum_{k = 1 - i}^0 |\phi_i| \\
&&\qquad\quad{}+ \|
\sqrt
{n} \overline{X}_n \|_p \Biggl(\| \overline{X}_n \|
_p +
n^{-1/2} \Biggl\| n^{-1/2}\sum_{k = 1}^{n - i} (X_k + X_{k + i})
\Biggr\|
_p\Biggr) \\
&&\qquad := A_n.
\end{eqnarray*}
Since $0 \leq i \leq d_n$, we obtain from Lemma~\ref{lemproparp} that
$A_n = \OO(n^{-1/2} d_n^{1/2})$, and hence by the Markov
inequality
%
%
\begin{eqnarray} \label{eqacomponentbound}
&&P\bigl(\bigl\|\sqrt{n} ({\widehat{\bolds\Phi}_{d_n}} - n^{-1}
{\mathbf
X}^T{\mathbf Y}) \bigr\|_{\infty}\geq\varepsilon\bigr)\nonumber\\[-8pt]\\[-8pt]
&&\qquad \leq\sum_{i =
1}^{d_n} P(|\Upsilon_i| \geq\varepsilon) = \OO
(\varepsilon
^{-p/2} n^{-p/4} d_n^{p/4 + 1}).\nonumber
\end{eqnarray}
Put $B_n = {\widehat{\bolds\Phi}_{d_n}} - n^{-1} {\mathbf X}^T{\mathbf
Y}$. Then
by adding and subtracting ${\bolds\Gamma}_{d_n}^{-1}$ we obtain
%
%
\begin{eqnarray}\label{eqacompoone}
P\bigl(\sqrt{n} \|{{\widehat{\bolds\Gamma}_{d_n} }}^{-1} B_n
\|
_{\infty}\geq\varepsilon\bigr)
&\leq& P\bigl(\sqrt{n}
\|
({{\widehat{\bolds\Gamma}_{d_n} }}^{-1} - {\bolds\Gamma
}_{d_n}^{-1}) B_n \|_{\infty}\geq\varepsilon/2 \bigr) \nonumber\\[-8pt]\\[-8pt]
&&{}
+P\bigl(\sqrt{n} \|{\bolds\Gamma}_{d_n}^{-1} B_n \|_{\infty
}\geq
\varepsilon/2 \bigr).\nonumber
\end{eqnarray}
In order to control the first expression, note that
\begin{eqnarray*}
P\bigl(\sqrt{n} \|({{\widehat{\bolds\Gamma}_{d_n} }}^{-1} -
{\bolds\Gamma}_{d_n}^{-1}) B_n \|_{\infty}\geq\varepsilon/2
\bigr) &\leq& P( \|{{\widehat{\bolds\Gamma}_{d_n} }}^{-1} - {\bolds
\Gamma
}_{d_n}^{-1}\|_{\infty} \varepsilon\geq\varepsilon/2 )\\
&&{} +
P(\|
B_n\|_{\infty} \geq\varepsilon),
\end{eqnarray*}
which by Lemma~\ref{lemmatrixinvapproxerror} and~(\ref{eqacomponentbound}) is of the magnitude $\OO(\varepsilon^{-p/2}
n^{-p/4} d_n^{p/4 + 1})$. Moreover, since $\|{\bolds\Gamma
}_{d_n}^{-1}\|_{\infty} < \infty$ by Lemma~\ref{lemboundmatrixnorms}, the bound in~(\ref{eqacomponentbound})
implies that for some $C > 0$
\[
P\bigl(\sqrt{n} \|{\bolds\Gamma}_{d_n}^{-1} B_n \|_{\infty
}\geq
\varepsilon/2 \bigr) \leq P\bigl(\sqrt{n} \|B_n \|_{\infty
}\geq
\varepsilon C^{-1} \bigr) =
\OO(\varepsilon^{-p/2} n^{-p/4} d_n^{p/4 + 1}),
\]
hence we conclude that
%
%
\begin{equation}\qquad
P\bigl(\bigl\|\sqrt{n}{{\widehat{\bolds\Gamma}_{d_n} }}^{-1}
({
\widehat{\bolds\Phi}_{d_n}} - n^{-1} {\mathbf X}^T{\mathbf Y})\bigr\|_{\infty}
\geq
\varepsilon\bigr) = \OO(\varepsilon^{-p/2} n^{-p/4} d_n^{p/4 +
1}).
\end{equation}

We will now treat the second part, which we rewrite as
\begin{eqnarray*}
&&n^{1/2}\bigl({{\widehat{\bolds\Gamma}_{d_n} }}^{-1} - n({\mathbf X}^T
{\mathbf X})^{-1} \bigr)\bigl( n^{-1} {\mathbf X}^T{\mathbf Y} - n^{-1}\E
({\mathbf X}^T{\mathbf Y})
\bigr)\\
&&\quad{}
+ n^{1/2}\bigl({{\widehat{\bolds\Gamma}_{d_n} }}^{-1} - n({\mathbf X}^T
{\mathbf X})^{-1} \bigr)n^{-1}\E({\mathbf X}^T{\mathbf Y}) \\
&&\qquad=: C_n + D_n.
\end{eqnarray*}

Due to Lemma~\ref{lemwuapprox} (requires an easy adaption), we have
%
%
\begin{equation}\label{eqwubound}
\|n^{-1/2} {\mathbf X}^T{\mathbf Y} - n^{-1/2}\E({\mathbf X}^T{\mathbf
Y}) \|
_{\infty}
= \OO_P( \log n ).
\end{equation}
Moreover, it holds that
\[
\sqrt{n}\bigl( {{\widehat{\bolds\Gamma}_{d_n} }}^{-1} - n({\mathbf X}^T
{\mathbf X})^{-1}\bigr) = {{\widehat{\bolds\Gamma}_{d_n} }}^{-1} \sqrt{n}
\bigl(n^{-1} ({\mathbf X}^T{\mathbf X}) - {{\widehat{\bolds\Gamma}_{d_n} }}
\bigr) n
({\mathbf
X}^T{\mathbf X})^{-1},
\]
and thus the sub-multiplicativity of the matrix norm \mbox{$\|\cdot\|_{\infty}$} implies
\begin{eqnarray*}
\|C_n\|_{\infty} &\leq&\|{{\widehat{\bolds\Gamma}_{d_n}
}}^{-1}\|_{\infty} \bigl\|\sqrt{n} \bigl(n^{-1} ({\mathbf X}^T{\mathbf X}) -
{{\widehat{\bolds\Gamma}_{d_n} }} \bigr)\bigr\|_{\infty} \|n
({\mathbf
X}^T{\mathbf X})^{-1}\|_{\infty} \\
&&{}\times\|n^{-1/2} {\mathbf X}^T{\mathbf Y} -
n^{-1/2}\E
({\mathbf X}^T{\mathbf Y})\|_{\infty}.
\end{eqnarray*}

Using~(\ref{eqwubound}) we thus obtain
\begin{eqnarray*}
&&
P(\|C_n\|_{\infty} \geq\varepsilon) \\
&&\qquad\leq
\OOO(1) +
P\bigl(\|{{\widehat{\bolds\Gamma}_{d_n} }}^{-1}\|_{\infty
}
\bigl\|\sqrt{n} \bigl(n^{-1} ({\mathbf X}^T{\mathbf X}) - {{\widehat{\bolds\Gamma
}_{d_n} }}
\bigr)\bigr\|_{\infty}  \\
&&\hspace*{-14.6pt}\hspace*{131pt}{}\times\|n ({\mathbf X}^T
{\mathbf X})^{-1}
\|_{\infty}\log n \geq\varepsilon\bigr).
\end{eqnarray*}

Put $\Delta_n = n^{-1} ({\mathbf X}^T{\mathbf X}) - {{\widehat{\bolds
\Gamma}_{d_n}}}$.
By adding and subtracting ${\bolds\Gamma}_{d_n}^{-1}$ and using Lem\-ma~\ref
{lemmatrixinvapproxerror} (see Remark~\ref{remmatrixinv}) and Lemma
\ref{lemboundmatrixnorms} we obtain
\begin{eqnarray*}
&&
P\bigl(\|{{\widehat{\bolds\Gamma}_{d_n} }}^{-1}\|_{\infty
}
\|\Delta_n \|_{\infty}  \|n ({\mathbf X}^T{\mathbf X})^{-1}\|
_{\infty
} \log n \geq\varepsilon\bigr) \\
&&\qquad\leq2 P( \|\Delta_n
\|
_{\infty} \log n \geq1 )  + P(\|{\bolds\Gamma
}_{d_n}^{-1} - \widehat{\bolds\Gamma}_{d_n}^{-1} \|_{\infty}
\geq
\varepsilon)\\
&&\qquad\quad{} + P\bigl(\|{\bolds\Gamma}_{d_n}^{-1} - n^{-1}
({\mathbf
X}^T{\mathbf X}) \|_{\infty} \geq\varepsilon\bigr).
\end{eqnarray*}
Choosing $\varepsilon= (\log n)^{-\chi_1}$, Lemma~\ref{lemmatrixinvapproxerror} and~(\ref{eqacomponentbound}) thus
yield the bound
%
%
\begin{equation}\qquad
P\bigl(\|{{\widehat{\bolds\Gamma}_{d_n} }}^{-1}\|_{\infty
}
\|\Delta_n \|_{\infty} \log n \geq(\log n)^{-\chi_1} \bigr)
= \OO
((\log n)^{\chi_1 p} n^{-p/2} d_n^{p/2}).
\end{equation}
Piecing everything together, the claim follows.
\end{pf}

We are now in the position to proof Lemma~\ref{lemapproxestimator}.
\begin{pf*}{Proof of Lemma~\ref{lemapproxestimator}}
We have that
\begin{eqnarray*}
&&
P\bigl(\| n^{1/2}(\widehat{\bolds\Theta} - {\bolds\Theta}) - n^{-1/2}
{{\bolds\Gamma}}^{-1} { \mathbf X}^T{\mathbf Z} \|_{\infty} \geq2
\varepsilon\bigr)
\\
&&\quad
\leq P\bigl(\| n^{1/2}(\widehat{\bolds\Theta} - \widetilde
{\bolds
\Theta})\|_{\infty} \geq\varepsilon\bigr) + P\bigl(\|
n^{1/2}(\widetilde{\bolds\Theta} - {\bolds\Theta}) - n^{-1/2} {{\bolds
\Gamma
}}^{-1} { \mathbf X}^T{\mathbf Z} \|_{\infty} \geq\varepsilon\bigr).
\end{eqnarray*}
Setting $\varepsilon= \log n^{-\chi_1}$, $\chi_1 > 2$, Proposition~\ref{propapproxestimatorreally} implies that
\[
\| n^{1/2}(\widehat{\bolds\Theta} - \widetilde{\bolds\Theta
})\|
_{\infty} = \OO_P(\log n^{-\chi_1} ).
\]
Moreover, the proof of Proposition~\ref{propapproxestimatorreally}
gives us
%
%
\begin{equation}\qquad
n^{1/2}(\widetilde{\bolds\Theta} - {\bolds\Theta}) - n^{-1/2} {{\bolds
\Gamma
}}^{-1} { \mathbf X}^T{\mathbf Z} = \bigl(n({\mathbf X}^T{\mathbf X})^{-1} -
{{\bolds\Gamma
}}^{-1}\bigr) n^{-1/2}{ \mathbf X}^T{\mathbf Z},
\end{equation}
and hence Remark~\ref{remmatrixinv} and Lemma~\ref{lemwuapprox}
imply that
%
%
\begin{eqnarray}
&&
\|n^{1/2}(\widetilde{\bolds\Theta} - {\bolds\Theta}) - n^{-1/2}
{{\bolds
\Gamma}}^{-1} { \mathbf X}^T{\mathbf Z} \|_{\infty} \nonumber\\
&&\qquad\leq\|n({\mathbf X}^T
{\mathbf X})^{-1} - {{\bolds\Gamma}}^{-1} \|_{\infty} \| n^{-1/2}{
\mathbf
X}^T{\mathbf Z} \|_{\infty} \\
&&\qquad= \OO_P(\log n^{-\chi_1 + 1} ),\nonumber
\end{eqnarray}
which completes the proof.
\end{pf*}

\section{Gaussian approximation}\label{secgaussianapprox}

In this section we obtain, under suitable assumptions, a~normal
approximation for the quantity
$n^{-1/2} \bolds\Gamma^{-1} \mathbf{X}^T \mathbf{Z}$, where we use the notation
introduced in Section~\ref{secyulewalkerdetails}. This entitles us to
obtain a~quantitative version of Theorem~\ref{thmyeconvrgextreme}
under stronger conditions. Let ${\mathbf V}_k = (X_{k - 1},\ldots, X_{k -
d_n})^T \varepsilon_k$, $k \in\N$. We have
\[
n^{-1/2} \bolds\Gamma^{-1} \mathbf{X}^T \mathbf{Z} = n^{-1/2}
{ \bolds\Gamma}^{-1} \sum_{k
= 1}^n {\mathbf V}_k = n^{-1/2} \sum_{k = 1}^n {\mathbf U}_k,
\]
where ${\mathbf V}_k = (V_{k}^{(1)},\ldots, V_{k}^{(d_n)})^T$,
${\mathbf
U}_k = (U_{k}^{(1)},\ldots, U_{k}^{(d_n)})^T$. Note that ${\mathbf
V}_k$ and~${\mathbf U}_k$ are both martingale sequences. In particular, it
holds that $\E({\mathbf V}_k ) = \E({\mathbf U}_k ) =
0$ and
%
%
\begin{equation}
\E({\mathbf V}_k {\mathbf V}_{k+h}^T ) = \cases{
\sigma^2 \Gamma_{d_n} , &\quad if $h = 0$,\cr
0_{d_n \times d_n}, &\quad if $h \neq0$,}
\end{equation}
since ${\bolds\varepsilon_k}$ is independent of $\{X_{k - i}\}_{i
\geq1}$. Throughout this section, we will always assume that $d_n =
\OO
(n)$.

The main theorem is formulated below.
%
%
\begin{theorem}\label{thmgaussianapprox}
Suppose that Assumption~\ref{assmainstrong} holds. If $d_n = \OO
(n^{\delta})$ with $\delta< 1/7$, then on a~possible larger
probability space, there exists a~$d_n$-dimen\-sional Gaussian random
vector ${\mathbf Z}$ with covariance matrix ${\bolds\Gamma}_{Z}$, such that
\[
P\Biggl(\Biggl\|{\mathbf Z} - \sum_{k = 1}^n {\mathbf U}_k \Biggr\|_{\infty}
\geq
v_n \Biggr) = \OO(n^{-\nu}),
\]
where $v_n = \sqrt{n} (\log n)^{-\chi_3}$, for arbitrary $\nu, \chi_3
\geq0$, and $\max\|n^{-1}{\bolds\Gamma}_{Z} - \sigma^2 {\bolds
\Gamma
}_{d_n} \| = \OOO(d_n^{-1}
)$.
\end{theorem}
%
%
\begin{rem}\label{remincreasedelta}
If one succeeds in establishing a~quantitative version of Lem\-ma~21 in
\cite{wuxiaobeats2011} with an appropriate error bound,
corresponding results to Theorem~\ref{thmgaussianapprox} with $0 <
\delta< 1$ should be possible. This, however, is beyond the scope of
the present paper.
\end{rem}

The proof of Theorem~\ref{thmgaussianapprox} partially follows
\cite{berkesgombayhi}, Theorem 4.1, and is based on a~series of lemmas. To
this end, we require some preliminary notation. For a~$d$-dimensional
vector ${\mathbf x} = (x_1,\ldots,x_d)$, we denote with $|{\mathbf
x}|_d =
(\sum_{i = 1}^n (x_i)^2)^{1/2}$ the usual Euclidean norm. The
following coupling inequality is due to Berthet and Mason
\cite{berthetmason}.
%
%
\begin{lem}[(Coupling inequality)]\label{lemcoupling}
Let $X_1,\ldots, X_N$ be independent, mean-zero random vectors in
$\R^n$, $n \geq1$, such that for some $B > 0$, $|X_i|_n \leq B$, $i =
1,\ldots, N$. If the probability space is rich enough, then for each
$\delta> 0$, one can define independent normally distributed mean-zero
random vectors $\xi_1,\ldots, \xi_N$ with $\xi_i$ and $X_i$ having the
same variance/covariance matrix for $i = 1,\ldots, N$, such that for
universal constants $C_1 > 0$ and $C_2 > 0$,
\[
P\Biggl\{\Biggl|\sum_{i = 1}^N (X_i - \xi_i)\Biggr|_n > \delta
\Biggr\}
\leq C_1 n^2 \exp\biggl( - \frac{C_2 \delta}{B n^2}\biggr).
\]
\end{lem}

The proof of Theorem~\ref{thmgaussianapprox} is based on a~blocking
argument, which in turn requires carefully truncated random variables. Put
\[
n^{-1/2} \bolds\Gamma^{-1} \mathbf{X}^T \mathbf{Z} = n^{-1/2}
{ \bolds\Gamma}^{-1} \sum_{k
= 1}^n {\mathbf V}_k = n^{-1/2} \sum_{k = 1}^n {\mathbf U}_k,
\]
where ${\mathbf U}_k = (U_{k}^{(1)},\ldots, U_{k}^{(d_n)})^T$. Note
that ${\mathbf V}_k$ and ${\mathbf U}_k$ are both martingale sequences.
%
%
\begin{lem}\label{lemwuapproxexact}
Suppose that Assumption~\ref{assmainstrong} holds. Then for $q \geq3$:
\begin{eqnarray*}
\mbox{\textup{(i)}\quad\hspace*{2pt}} P\Biggl(\Biggl\|n^{-1/2}\sum_{k = 1}^n {\mathbf U}_{k} \Biggr\|
_{\infty
} \geq\sqrt{q \log n} \Biggr) &=& \OO(n^{-\nu}),\\
\mbox{\textup{(ii)}\quad} P\bigl(\sqrt{n}\|\widehat{\bolds\Phi}_{d_n} - {\bolds\Phi
}_{d_n} \|_{\infty} \geq\sqrt{q \log n}\bigr) &=& \OO
(n^{-\nu})
\end{eqnarray*}
for arbitrary $\nu\geq0$.
\end{lem}
\begin{pf}
We first show (i). By Lemma 1 in~\cite{wuasymptoticbernoulli} we have
\begin{eqnarray*}
P\Biggl(\Biggl\|n^{-1/2}\sum_{k = 1}^n {\mathbf U}_{k} \Biggr\|_{\infty}
\geq
\sqrt{q \log n} \Biggr) &\leq&\sum_{h = 1}^{d_n} P\Biggl(
\Biggl|n^{-1/2}\sum
_{k = 1}^n {U}_{k}^{(h)} \Biggr| \geq\sqrt{q \log n} \Biggr) \\
&=& \OO
(d_n n^{-\nu})
\end{eqnarray*}
for arbitrary $\nu\geq0$, hence the claim. Part (ii) can be shown in
the same way, using Theorem 3 in~\cite{wuasymptoticbernoulli}
instead of Lemma 1.
\end{pf}
%
%
\begin{lem}\label{lemapproxandbound}
If Assumption~\ref{assmainstrong} is valid, then there exists a~sequence of random vectors ${\mathbf U}_k^* = (U_{k}^{(1,*)},\ldots,
U_{k}^{(d_n,*)})^T$ with $\E({\mathbf U}_k^* ) = 0$ and
the same covariance structure as ${\mathbf U}_k$, such that ${\mathbf U}_k^*$
is a~$d_n$-dependent sequence,\break ${\max_{1 \leq k \leq n}}
|U_{k}^{(h,*)} | = \OO(b_n^2)$, $1 \leq h \leq d_n$, and
\[
P\Biggl(n^{-1/2}\Biggl\|\sum_{k = 1}^n {\mathbf U}_k - \sum_{k = 1}^n
{\mathbf
U}_{k}^{*} \Biggr\|_{\infty} \geq v_n \Biggr) = \OO(n^{-\nu
}),
\]
where $v_n = \sqrt{n} (\log n)^{-\chi_3}$ for arbitrary $\nu, \chi_3
\geq0$.
\end{lem}
\begin{pf}
Put
%
%
\begin{equation}
\varepsilon_{k,b_n} = \varepsilon_k\ind_{|\varepsilon_k| \leq b_n} - \E
\bigl(\varepsilon_k \ind_{|\varepsilon_k| \leq b_n} \bigr)
\end{equation}
and let
\begin{eqnarray*}
U_{{k,b_n}}^{(h)} &=& U_{k}^{(h)}\ind_{{\max_{|l|\leq n}}|\varepsilon_l|
\leq
b_n} - \E\bigl(U_{k}^{(h)}\ind_{{\max_{|l|\leq n}}|\varepsilon_l| \leq b_n}
\bigr) \\
&=& \varepsilon_{k,b_n} \Biggl(\sum_{j = 1}^{d_n} \gamma_{h,j}^{(*)}
\sum_{i = 0}^{\infty} \alpha_i \varepsilon_{k - j - i,b_n} \Biggr).
\end{eqnarray*}
Denote with ${\mathbf U}_{k,b_n} = (U_{{k,b_n}}^{(1)},\ldots,
U_{{k,b_n}}^{(d_n)})^T$; then
\begin{eqnarray*}
&&P\Biggl(\Biggl\|\sum_{k = 1}^n {\mathbf U}_k - \sum_{k = 1}^n {\mathbf U}_{k,b_n}
\Biggr\|_{\infty} \geq v_n \Biggr) \\
&&\qquad\leq P\Bigl({\max_{|l| \leq
n}}|\varepsilon
_l| > b_n \Bigr) + P\bigl(\bigl|\sqrt{n} \E\bigl({\mathbf U}_{k,b_n}^{(h)}
\bigr) \bigr|
\geq(\log n)^{-\chi_3} \bigr).
\end{eqnarray*}

Since $\E({\mathbf U}_{k}^{(h)} ) = 0$, an application of the
Cauchy--Schwarz inequality yields
\[
\bigl|\sqrt{n} \E\bigl({\mathbf U}_{k,b_n}^{(h)} \bigr) \bigr| \leq
\sqrt
{n} \bigl\|{\mathbf U}_{k,b_n}^{(h)} \bigr\|_2 \bigl\|\ind_{{\max
_{|l|\leq
n}}|\varepsilon_l| > b_n} \bigr\|_2 = C \sqrt{n P\Bigl({\max_{|l|\leq
n}}|\varepsilon_l| > b_n \Bigr)},
\]
which by Assumption~\ref{assmainstrong} is of the magnitude $\OO
(n^{-\nu})$, for arbitrary $\nu\geq0$. Hence we conclude
%
%
\begin{equation}\label{eqboundok}
P\Biggl(\Biggl\|\sum_{k = 1}^n {\mathbf U}_k - \sum_{k = 1}^n {\mathbf U}_{k,b_n}
\Biggr\|_{\infty} \geq v_n \Biggr) = \OO(n^{-\nu}).
\end{equation}
Put
${\mathbf U}_{k,b_n}^{(d_n)} = (U_{{k,b_n}}^{(1,d_n)},\ldots,
U_{{k,b_n}}^{(d_n,d_n)})^T$. Then
\[
{\mathbf U}_{k,b_n}^{(d_n)} = \varepsilon_{k,b_n} \Biggl(\sum_{j = 1}^{d_n}
\gamma_{h,j}^{(*)} \sum_{i = 0}^{d_n} \alpha_i \varepsilon_{k - j - i,b_n}
\Biggr).
\]

By Lemma~\ref{lemwuapproxexact} (remains valid) we have that
\begin{eqnarray*}
&&
P\Biggl(\Biggl\|\sum_{k = 1}^n {\mathbf U}_{k,b_n} - {\mathbf U}_{k,b_n}^{(d_n)}
\Biggr\|_{\infty} \geq v_n \Biggr) \\
&&\qquad\leq
\sum_{h = 0}^{d_n} P
\Biggl(\Psi
(d_n)^{-1/2} \Biggl|n^{-1/2}\sum_{k = 1}^n {U}_{k,b_n}^{(h)} - {
U}_{k,b_n}^{(h,d_n)} \Biggr| \geq\Psi(d_n)^{-1/2} (\log n)^{-\chi_3}
\Biggr) \\
&&\qquad= \OO(n^{-\nu})
\end{eqnarray*}
for arbitrary $\nu\geq0$. Let $\{\varepsilon_k^{(h,*)}\}_{k
\in\Z}$, $1 \leq h \leq\mathfrak{d}$, be an array of mutually
independent random variables, where $\varepsilon_k^{(h,*)}$ is an
independent copy of $\varepsilon_{k,b_n}$ for each $h$. Then we can define
the random vectors
\[
{U}_{k}^{(h,*)} = {U}_{k,b_n}^{(h,d_n,*)} = \varepsilon_{k,b_n}
\Biggl(\sum
_{j = 1}^{d_n} \gamma_{h,j}^{(*)} \Biggl[\sum_{i = 0}^{d_n} \alpha_i
\varepsilon_{k - j - i,b_n} + \sum_{i = d_n + 1}^{\infty} \alpha_i
\varepsilon
_{k - j - i}^{(h,*)} \Biggr]\Biggr).
\]
Note\vspace*{1pt} that due to the structure of
${U}_{k,b_n}^{(h,d_n,*)}$ it is clear that one may repeat all the
previous arguments to derive the bound
%
%
\begin{equation}
n^{-1/2} \Biggl\|\sum_{k = 1}^n {\mathbf U}_k - \sum_{k = 1}^n
{U}_{k,b_n}^{(h,d_n,*)} \Biggr\|_{\infty} = \OO_P(n^{-\nu
}).
\end{equation}
Let $\sigma_n^{*} = \Var(\varepsilon_{k,b_n})$. Since
$\sigma
_n^{*} > 0$ for large enough $n$, the Cauchy--Schwarz inequality and
Assumption~\ref{assmain} imply
\begin{eqnarray*}
\bigl|\sqrt{\sigma_n^{*}} - \sqrt{\sigma^{2}} \bigr|_1 &\leq& C
|\sigma_n^{*} - \sigma^{2} |_1 = C \bigl\|\varepsilon_k^2 \ind
_{|\varepsilon_k| \geq b_n} \bigr\|_1 \\
&\leq& C \|\varepsilon_k^2
\|_2
\sqrt{P(|\varepsilon_k| \geq b_n )} = \OO(n^{-\nu
}).
\end{eqnarray*}

Then we obtain from the above and Lemma~\ref{lemwuapproxexact}
(remains valid)
%
%
\begin{equation}
n^{-1/2}\Biggl\| (1 - \sigma^2/ \sigma_n^{*})\sum_{k = 1}^n
{U}_{k,b_n}^{(h,d_n,*)} \Biggr\|_{\infty} = \OO_P(n^{-\nu
}).
\end{equation}
Put\vspace*{1pt} ${\mathbf U}_k^* = ({U}_{k}^{(1,*)},\ldots,{U}_{k}^{(d_n,*)}
)^T$. Then it is clear that
${\max_{1 \leq k \leq n}}|U_{k}^{(h,*)} |_d = \OO
(b_n^2)$, $1 \leq h \leq d_n$,
and piecing everything together, the claim follows.
\end{pf}

We will now construct an approximation for the random vector ${\mathbf
U}_{k}^*$.
To this end, we first divide the set of integers
$\{1,2,\ldots\}$ into consecutive blocks $H_1, J_1$, $H_2,J_2,\ldots.$ The
blocks are defined by recursion. Fix $\delta^* > \delta> 0$, and put
$m_n = \lfloor n^{\delta^*} \rfloor$. If the largest
element of $J_{i-1}$ is $k_{i-1}$, then $H_i = \{k_{i-1} +
1,\ldots,k_{i-1} + m_n\}$ and $J_i = \{k_{i-1} + m_n +
1,\ldots, d_n\}$. Let $|\cdot|$ denote the cardinality of a~set. It
follows from the definition of $H_i$, $J_i$ that $|H_i| = m_n$
and $|J_i| = d_n$. Note that the\vadjust{\goodbreak} total number of blocks is
approximately $n/m_m = n^{1 - \delta^*}$. Let $\mathcal{I} \subset\{
0,1,\ldots,d_n\}$ be a~subset with $|\mathcal{I}| = \mathfrak{d}$, with
$\mathfrak{d} = \OO(n^{\lambda})$, $\lambda> 0$, and denote
with $\sigma^2 {\bolds\Gamma}_{\mathcal{I}}$ the sub-covariance matrix of
${\mathbf U}_{k}^*$ restricted to the subset $\mathcal{I}$.
%
%
\begin{lem}\label{lemgaussianapprox}
If Assumption~\ref{assmainstrong} is valid and $5\lambda+ 2 \delta^*
< 1$, then on a~possible larger probability space there exists a~$\mathfrak{d}$-dimensional Gaussian random vector ${\mathbf Z}$ with
covariance matrix $n {\bolds\Gamma}_{Z,\mathcal{I}}$, such that
\[
P\Biggl(\max_{h \in\mathcal{I}}\Biggl|{\mathbf Z} - \sum_{k = 1}^n
{\mathbf
U}_k^* \Biggr| \geq v_n \Biggr) = \OO(\exp(-n^{\varepsilon
})),\qquad
\varepsilon> 0,
\]
where $v_n = \sqrt{n} (\log n)^{-\chi_3}$, for arbitrary $\chi_3
\geq
0$, and ${\max}\|{\bolds\Gamma}_{Z,\mathcal{I}} - \sigma^2 {\bolds
\Gamma
}_{\mathcal{I}} \| = \OO(m_n^{-1} )$.
\end{lem}
\begin{pf}
For $h \in\mathcal{I}$, let
\[
\xi_k^{(h)} = \sum_{i \in H_k} U_{i}^{(h,*)} \quad\mbox{and}\quad
\eta
_k^{(h)} = \sum_{i \in J_k} U_{i}^{(h,*)}
\]
and define the vectors
\[
\xiv_k = \bigl(\ldots, \xi_k^{(h)},\ldots\bigr)^{T} ,\qquad h \in\mathcal{I}
\mbox{,\quad and}\quad \etav_k = \bigl(\ldots,
\eta_k^{(h)},\ldots\bigr)^{T},\qquad
h\in\mathcal{I}.
\]
Note that per construction, we have that $\{\xiv_k \}_{k
\in
\N}$ is a~sequence of independent random vectors with $|\xiv
_k|_{\mathfrak{d}} = \OO(\sqrt{\mathfrak{d}} m_n b_n^2
)$. By
Lemma~\ref{lemcoupling}, we can define a~sequence of independent
normal random vectors $\xiv_k^* = (\ldots,\xi_k^{(h,*)},\ldots
)^{T}$, $h \in\mathcal{I}$, such that for $x > 0$
\begin{eqnarray*}
P\Biggl(\max_{1 \leq h \leq\mathfrak{d}}
\Biggl|\sum_{j = 1}^{n/m_n} \bigl(\xi_j^{(h)} - \xi_j^{(h,*)}\bigr)
\Biggr| \geq x \Biggr) & \leq & \sum_{h = 1}^{\mathfrak{d}}
P\Biggl( \Biggl|\sum_{j = 1}^{n/m_n} \bigl(\xi_j^{(h)} -
\xi_j^{(h,*)}\bigr) \Biggr| \geq x \Biggr)
\\
&\leq& \sum_{h = 1}^{\mathfrak{d}} P\Biggl(\Biggl|\sum_{j = 1}^{n/m_n}
({\xiv_j} - {\xiv_j^*}) \Biggr|_{\mathfrak{d}} \geq x
\Biggr)
\\
&\leq& C \mathfrak{d}^2 \exp\biggl(-\frac{x}{\mathfrak{d}^{5/2}
m_n b_n^2}
\biggr).
\end{eqnarray*}

We thus obtain
%
%
\begin{equation}
P\Biggl(\max_{1 \leq h \leq\mathfrak{d}} \Biggl|\sum_{j = 1}^{n/m_n}
\bigl(\xi_j^{(h)} - \xi_j^{(h,*)}\bigr) \Biggr| \geq v_n \Biggr)
= \OO
(\exp(-n^{\varepsilon})),
\end{equation}
and similar arguments show that there exists a~sequence of independent normal
random vectors $\etav_k^* = (\ldots,\eta_k^{(h,*)},\ldots)^{T}$,
such that
\[
P\Biggl(\max_{1 \leq h \leq\mathfrak{d}}
\Biggl|\sum_{j = 1}^{n/m_n} \bigl(\eta_j^{(h)} - \eta
_j^{(h,*)}\bigr)
\Biggr| \geq v_n \Biggr)= \OO(\exp(-n^{\varepsilon})).
\]
Lemma~\ref{lemapproxandbound} yields that $\Var(
\eta_j^{(h,*)} ) =\OO(d_n)$ for all $j \leq m_n$, $1
\leq h \leq\mathfrak{d}$. Hence by known
properties of the tails of a~normal c.d.f., we obtain that
%
%
\begin{eqnarray}\label{eqmaxboundcomp}
P\Biggl(\max_{1 \leq h \leq\mathfrak{d}}
\Biggl|\sum_{j = 1}^{n/m_n} \eta_j^{(h),*}\Biggr| \geq v_n \Biggr)
&\leq& \sum_{h = 1}^{\mathfrak{d}} P\Biggl( \Biggl|\sum_{j = 1}^{n/m_n}
\eta_j^{(h),*}\Biggr| \geq v_n \Biggr) \nonumber\\
&\leq& {\mathfrak{d}}
P
\bigl(|Z| \geq C \sqrt{d_n/m_n} (\log n)^{-\chi_3}\bigr) \\
&=& \OO
(\exp
(-n^{\varepsilon}))\nonumber
\end{eqnarray}
for some $\varepsilon> 0$.
This yields
%
%
\begin{equation}
P\Biggl(\max_{1 \leq h \leq\mathfrak{d}} \Biggl|\sum_{j = 1}^{n/m_n}
\bigl(\eta_j^{(h)} + \xi_j^{(h)} - \xi_j^{(h,*)} \bigr) \Biggr|
\geq
v_n \Biggr)= \OO(\exp(-n^{\varepsilon})).
\end{equation}
Let $\etav_k^{**} = (\ldots,\eta_k^{(h,**)},\ldots)^{T}$ $h
\in
\mathcal{I}$ be a~copy of
$\etav_k^{*}$ such that $\etav_i^{**}$ and $\xiv_j^{*}$
are independent for $i \neq j$. By the very construction of $\xiv_k,
\etav_k$, it is not hard to show that
\begin{eqnarray*}
&&\max_{i,j \in\mathcal{I}}\Biggl|\Cov\Biggl(\sum_{k = 1}^{n/m_n}
\eta
_k^{(i)} + \xi_k^{(i)}, \sum_{k = 1}^{n/m_n} \eta_k^{(j)} + \xi_k^{(j)}
\Biggr) \\
&&\hspace*{10pt}\quad{}- \Cov\Biggl(\sum_{k = 1}^{n/m_n} {\xi}_{k}^{(i,*)} +
\eta
_k^{(i,**)}, \sum_{k = 1}^{n/m_n} {\xi}_{k}^{(j,*)} + \eta_k^{(j,**)}
\Biggr) \Biggr| \\
&&\qquad= \OO(n/m_n),
\end{eqnarray*}
which clearly implies ${\max}\|{\bolds\Gamma}_{Z,\mathcal{I}} -
\sigma
^2{\bolds\Gamma}_{\mathcal{I}} \| = \OO(m_n^{-1} )$.
Hence, by enlarging the probability space if necessary and arguing
similarly as in~(\ref{eqmaxboundcomp}),
we have that
\[
P\Biggl(\max_{1 \leq h \leq\mathfrak{d}}
\Biggl|\sum_{j = 1}^{n/m_n} \bigl(\xi_j^{(h)} + \eta_j^{(h)} - \xi
_j^{(h,*)} - \eta_j^{(h,**)}\bigr) \Biggr| \geq v_n \Biggr) = \OO
(\exp(-n^{\varepsilon})).
\]
Finally, we obtain from the above
\[
P\Biggl(\max_{h \in\mathcal{I}} \Biggl|\sum_{k = 1}^n {\mathbf U}_k^* -
\sum_{j = 1}^{n/m_n} ( \xiv_j^{*} - \etav_j^{**}) \Biggr|
\geq v_n \Biggr) = \OO(\exp(-n^{\varepsilon})),
\]
which completes the proof.
\end{pf}
\begin{pf*}{Proof of Theorem~\ref{thmgaussianapprox}}
By Lemma~\ref{lemapproxandbound} it suffices to establish the claim
for $\{{\mathbf U}_k^*\}_{1 \leq k \leq n}$. This, however, is
provided by Lemma~\ref{lemgaussianapprox}.
\end{pf*}

\section*{Acknowledgments}

The author thanks the Associate Editor and the anon\-ymous referee for
their valuable remarks and comments that considerably helped to improve
the results of this paper.\vadjust{\goodbreak}


%

\printaddresses

\end{document}